\documentclass[a4paper,reqno,12pt]{amsart}
\usepackage{rotating,pdflscape,color,amssymb,comment,subfigure,epstopdf,psfrag,amsmath,bbm,eufrak,graphicx}

\usepackage{mca_prefs}
\author[Raj Rao Nadakuditi]{Raj Rao Nadakuditi}\address{Raj Rao Nadakuditi, Department of Electrical Engineering and Computer Science, University of Michigan, 1301 Beal Avenue, Ann Arbor, MI 48109. USA.}
\email{rajnrao@eecs.umich.edu}
\urladdr{http://www.eecs.umich.edu/\~\/rajnrao/}

\title[Informative component analysis]{When are the most informative components for inference also the principal components? }

\thanks{This work was supported by the ARO MURI W911NF-11-1-0391 grant and an AFOSR Young Investigator Award FA9550-12-1-0266. The author thanks Florent Benaych-Georges for many inspiring conversations that led to this work and Iain Johnstone for suggesting Gaussian mixture models as an example where middle components might be informative.  }

\date{\today}

\begin{document}
\maketitle
\begin{abstract}Which components of the singular value decomposition of a signal-plus-noise data matrix are most informative for the inferential task of detecting or estimating an embedded low-rank signal matrix? Principal component analysis ascribes greater importance to the components that capture the greatest variation, \textit{i.e.}, the singular vectors associated with the largest singular values.  This choice is often justified by invoking the Eckart-Young theorem even though that work addresses the problem of how to best \textit{represent} a signal-plus-noise matrix using a low-rank approximation and \textit{not} how to best \textit{infer} the underlying low-rank signal component.

Here we take a first-principles approach in which we start with a signal-plus-noise data matrix and show how the spectrum of the noise-only component governs whether the principal or the middle components of the singular value decomposition of the data matrix will be the informative components for inference.

Simply put, if the noise spectrum is supported on a connected interval, in a sense we make precise, then the use of the  principal components is justified. When the noise spectrum is supported on multiple intervals, then the middle components might be more informative than the principal components.

The end result is a proper justification of the use of principal components in the oft considered setting where the noise matrix is i.i.d. Gaussian. An additional consequence of our study is  the identification of scenarios, generically involving heterogeneous noise models such as mixtures of Gaussians, where the middle components might be more informative than the principal components so that they may be exploited to extract additional processing gain. In these settings, our results show how the blind use of principal components can lead to suboptimal or even faulty inference because of phase transitions that separate a regime where the principal components are informative from a regime where they are uninformative. We illustrate our findings using numerical simulations and a real-world example.
\end{abstract}


\section{Introduction}
Consider a signal-plus-noise data matrix modeled as
\begin{equation}\label{eq:sigplusnoise mat}
\wtX = \sum_{i=1}^{r} \theta_{i} u_{i} v_{i}^{H} + X,
\end{equation}
where $X$ denotes the $n \times m$ noise-only matrix and $S = \sum_{i=1}^{r} \theta_{i} u_{i} v_{i}^{H}$ is the rank-$r$ signal matrix.  Relative to this model, the detection and estimation tasks in signal processing and data analysis deal with inferring the presence of and estimating the rank $r$ matrix $S$ given $\wtX$.

Principal component analysis plays an important role in the setting where $r \ll \min (m,n)$ as described succinctly by Joliffe \cite[Ch1., pp.1]{jolliffe2005principal}:
\begin{quote}
The central idea of principal component analysis (PCA) is to reduce the dimensionality of a data set .... while \emph{retaining as much as possible of the variation}\footnote{Emphasis added.} present in the data set. The ... \emph{first few retain most of the variation}\footnotemark[\value{footnote}] present  in all of the original variables.
\end{quote}

The \emph{first few} principal components alluded to here refer to the first few singular vectors associated with the largest singular values of $\wtX$.
Working with the hypothesis that the directions of greatest variation of the data set must reflect (or correlate with) the signal content and equipped with the singular value decomposition (SVD) as a technique for computing these directions, we can tackle the detection problem in the following manner.

We start off by computing the SVD of $\wtX$ and plot the singular values $\{\widetilde{\sigma}_{i}\}_{i=1}^{n}$ in non-increasing order. We then estimate the rank of the latent signal matrix $S$ based on the rule:
\begin{equation}\label{eq:khat}
 \widehat{r} = \{\textrm{First} \,i \textrm{ such that }  \textrm{gap}(i):= \widetilde{\sigma}_{i} - \sigma_{\sf null.} < \textrm{threshold} \} -1,
 \end{equation}
where $\sigma_{\sf null}$ is the largest singular value of the noise-only matrix $X$ which is assumed (in the simplest setting) to be known. This rule, and other modifications thereof, yields an estimate $\widehat{r}$ for the rank of the latent signal matrix; when $\widehat{r} >0$ we have detected a signal matrix; see for example \cite[Section 14.5]{friedman2001elements} or \cite[Section 6.1.3]{jolliffe2005principal} for classical approaches and \cite{johnstone2001distribution,baik2005phase,baik2006eigenvalues,johnstone2006high,el2007tracy,paul2007asymptotics,nadakuditi2008sample,onatski2010determining,kritchman2008determining,kritchman2009non,nadler2010nonparametric} for recent random matrix-theoretic approaches.

\begin{figure}[t]
\centering
\subfigure[The singular value spectrum.]
{
\includegraphics[width=5.25in]{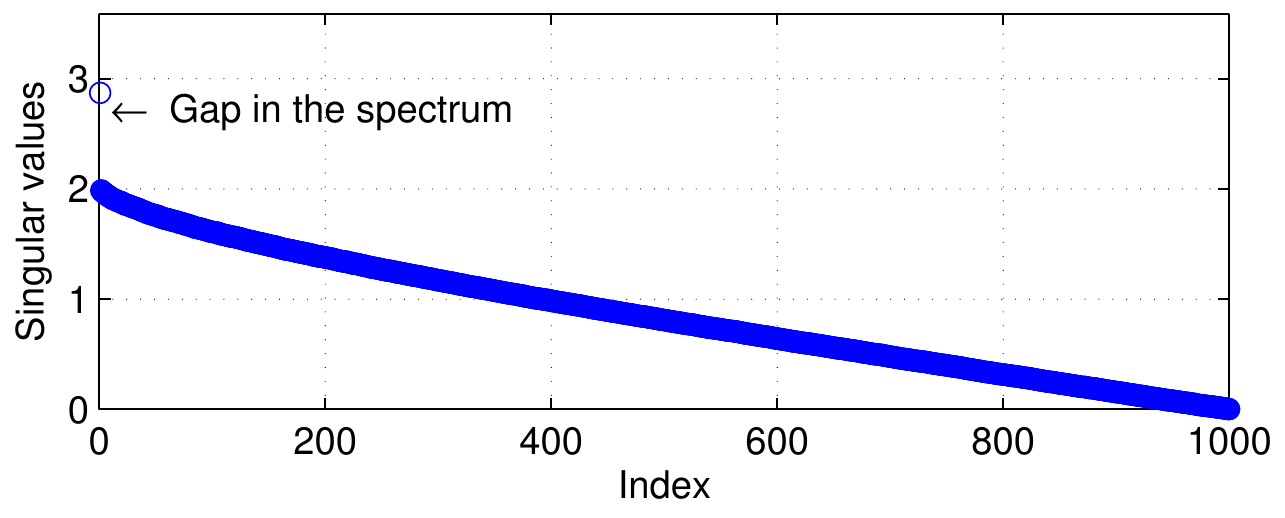}
\label{fig:pca sv}
}\\
\subfigure[PCA example.]
{
\includegraphics[width=5.25in]{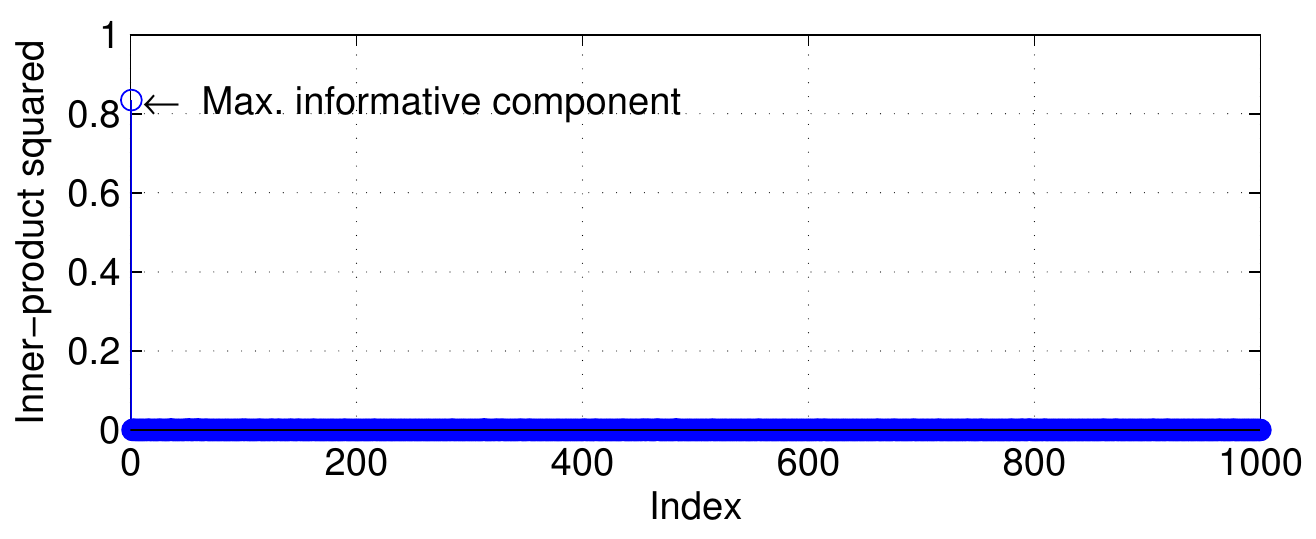}
\label{fig:pca ip}
}
\label{fig:typical pca}
\caption{The singular value spectrum of the signal-plus-noise data matrix exhibits a continuous-looking portion that may be associated with ``noise'' singular values and a single separated singular value that may be interpreted as evidence of a rank-one ``signal'' matrix buried in the data matrix.}
\end{figure}

The estimation problem is similarly tackled by computing the truncated SVD of $\wtX$ that employs the $\widehat{r}$ (leading or) principal components. This yields a rank $\widehat{r}$ estimate of the low-rank signal matrix given by
\begin{equation}\label{eq:Shat}
\widehat{S} = \sum_{i=1}^{\widehat{r}} \widetilde{\sigma}_{i} \widetilde{u}_{i} \widetilde{v}_{i}^{H}.
\end{equation}
Does the principal component approach to detection and estimation work? Figure \ref{fig:pca sv} plots the singular values of a $n \times m$ signal-plus-noise data matrix modeled as $\wtX = 2 uv^{H} + X$, where the noise-only matrix $X$ has i.i.d. mean zero, variance $1/m$ Gaussian entries and the signal matrix $S = 2 uv^{H}$ has rank one. This example, where $n = m = 1000$, illustrates a setting  where the gap heuristic in (\ref{eq:khat})  for signal-matrix detection ``works''  subject to a specification of the gap size threshold.

Figure \ref{fig:pca ip} plots the $n$ inner-products $\{|\langle \widetilde{u}_{i},u \rangle|^{2}|\}_{i=1}^{n}$, where $\{\widetilde{u}_{i}\}_{i=1}^{n}$ are the left singular vectors of $\wtX$. The quantities  $\{|\langle \widetilde{u}_{i},u \rangle|^{2}|\}_{i=1}^{n}$ (and $\{|\langle \widetilde{v}_{i},v \rangle|^{2}|\}_{i=1}^{m}$) are measures of informativeness of the singular vectors of $\wtX$ with respect to the singular vectors of the latent signal matrix. Clearly, the principal left (also, the right - not plotted here) singular vector is the most informative component  and employing it in an estimate of the signal matrix as in (\ref{eq:Shat}) is judicious.

Extending the notion of informativeness further, we might define ``informative components'' as components of the SVD of the data matrix $\wtX$ that are most correlated with the embedded low-rank signal matrix and which consequently best (in a manner to be made precise later) facilitate the detection and estimation tasks described earlier.

For the example in Figure \ref{fig:pca ip}, the {principal component} \textit{is} the most {informative component}. In other words, the principal component which captures the greatest \textit{variation} in the data is also the component most correlated with the underlying signal matrix. A natural question arises:
\begin{center}
Are the most \textit{informative components} necessarily the \textit{principal components}?
\end{center}

\begin{figure}[t]
\centering
\subfigure[MCA example.]
{
\includegraphics[width=5.25in]{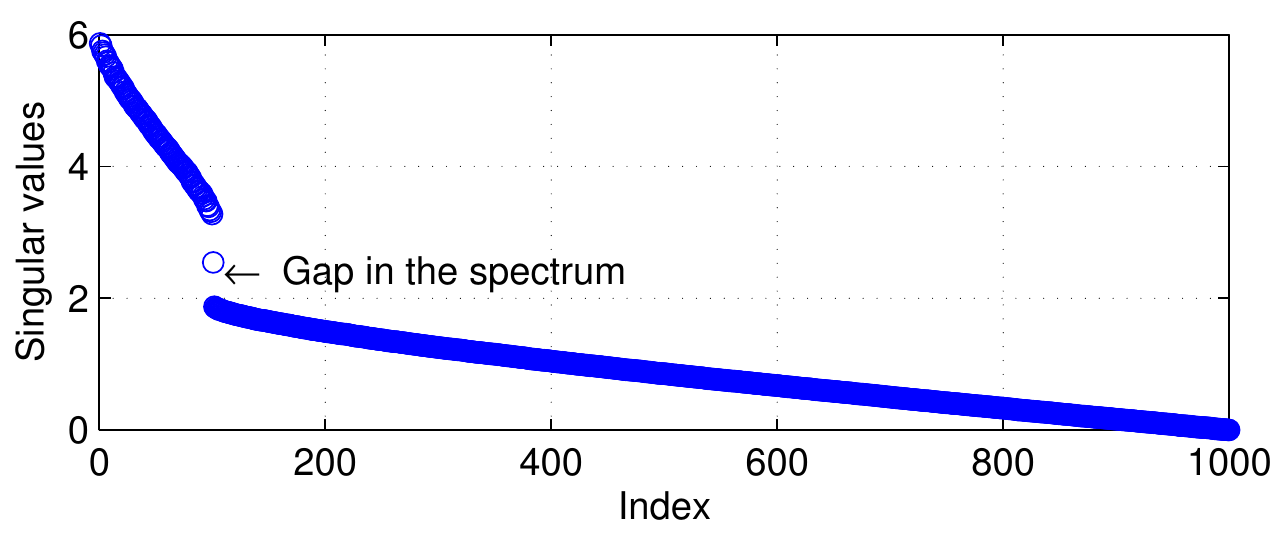}
\label{fig:mca sv}

}
\subfigure[MCA example.]
{
\includegraphics[width=5.25in]{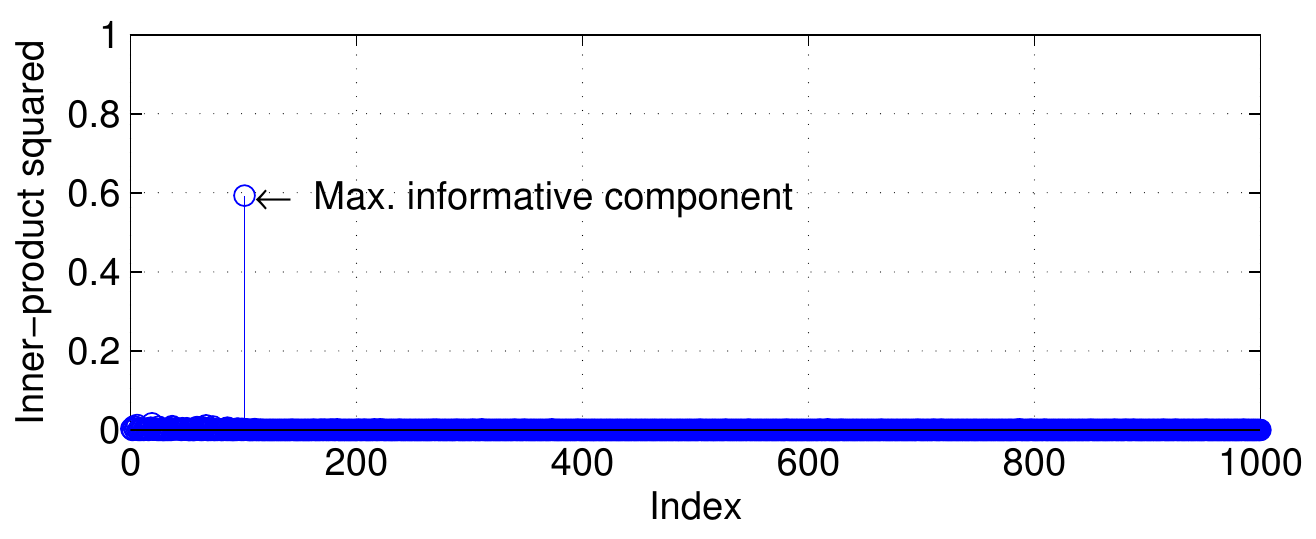}
\label{fig:mca ip}
}
\caption{The singular value spectrum of the signal-plus-noise data matrix exhibits two continuous-looking portions that may be associated with ``noise'' singular values and a single separated singular value that may be interpreted as evidence of a rank-one ``signal'' matrix buried in the data matrix. Note that in contrast to Figure \ref{fig:typical pca}, the principal component is \textit{not} the most informative component.}
\label{fig:typical mca}
\end{figure}

Figure \ref{fig:typical mca} constitutes a counter-example. Figure \ref{fig:mca sv} plots the singular values of a signal-plus-noise data matrix modeled as $\wtX = 2 uv^{H} + X$, where the noise-only matrix $X$ is a mixture of two multivariate Gaussians with different variances that produces a spectrum that is supported on two disconnected intervals. The MATLAB code used to generate $\wtX$ is listed below so the reader may reproduce Figure \ref{fig:typical mca}:

\begin{verbatim}
n = 1000; m = n;
Sigma = diag([20*ones(n/10,1);ones(n-n/10,1)],0); % temporal covariance
G = randn(n,m)/sqrt(m)*sqrtm(Sigma);
u = randn(n,1); u = u/norm(u); v = randn(m,1)/sqrt(m);
Xtil = 2*u*v' + G;
\end{verbatim}

The presence of a signal matrix is reflected in the single singular value that separates from the continuous looking portions of the spectrum - unlike Figure \ref{fig:pca sv}, it is in the middle \textit{i.e.} \textbf{not} associated with the principal component that captures the greatest variation. The rule in (\ref{eq:khat}) would return $\widehat{r} = 0$ here and we would fail to detect the underlying signal matrix.

Figure \ref{fig:mca ip} plots the inner-product $\{|\langle \widetilde{u}_{i},u \rangle|^{2}|\}_{i=1}^{n}$, where $\{\widetilde{u}_{i}\}_{i=1}^{n}$ are the left singular vectors of $\wtX$. The quantities  $\{|\langle \widetilde{u}_{i},u \rangle|^{2}|\}_{i=1}^{n}$ (and $\{|\langle \widetilde{v}_{i},v \rangle|^{2}|\}_{i=1}^{m}$) are measures of informativeness of the singular vectors of $\wtX$ with respect to the singular vectors of the latent signal matrix. Clearly, the principal left (also, the right - not plotted here) singular vector is \textbf{not} the most informative component; the middle component is. Employing the principal component in an estimate of the signal matrix as in (\ref{eq:Shat}) would not be as judicious as using the most informative component, which is the \textit{middle component} here.

The preceding examples support our assertion that the principal components  are not necessarily the most informative components and that middle components might sometimes be more important. The examples also hint at the role played by the spectrum of the noise-only matrix $X$ in determining the relative informativeness of the components.

An additional remark is in order. The Eckart-Young-Mirsky (EYM) theorem \cite{eckart1936approximation,mirsky1960symmetric} states that for any unitarily invariant norm, the optimal rank $\widehat{r}$ approximation to $\widehat{X}_{n}$ is given by (\ref{eq:Shat}).  This is a statement about optimal \textit{representation} of the signal-plus-noise matrix. It is \textbf{not} a statement about inference on the underlying low-rank signal matrix. Thus there is no contradiction between our results and the content of the EYM theorem.

\subsection{Motivation and summary of findings}
 This work is motivated by the ubiquity of principal component analysis (PCA) in data analysis and signal processing and the associated importance assigned by practitioners to the leading singular values and vectors of the data matrix.

In emerging applications, such as the collaborative learning, graph mining or bioinformatics where the data matrix is large, it is infeasible to compute the entire singular value decomposition. There are, however, efficient techniques for computing the leading singular vectors of a matrix that employ iterative techniques such as the Arnoldi or Lanczos iteration \cite{brand2006fast} and the family of Krylov subspace methods or using randomized techniques as in \cite{frieze2004fast,drineas2005nystrom,deshpande2006adaptive,halko2011finding}.

In these `big data' applications, researchers often invoke PCA as justification for the computation of a small number of leading singular vectors of the data matrix. Arguably, what a practitioner who uses these principal components as a starting point in an inferential detection, estimation or classification procedure is really after are the informative components. As we have already seen, the informative components need not be the principal components and may even be the middle components.

In the latter scenario, computation of the leading singular vectors, regardless of computational considerations or choice of algorithm, might lead to faulty inference and lead a non-specialist down a road to a flawed conclusion that they may present as supported by standard PCA derived data analysis.  The situation is particularly perilous in biomedical applications involving high-dimensional data sets where one cannot exclude or reason about most informative components by visual inspection.  \footnote{\url{http://www.nytimes.com/2011/07/19/health/19gene.html?pagewanted=all}}\textsuperscript{,}\footnote{\url{http://www.nytimes.com/2011/07/08/health/research/08genes.html?_r=2&hp}}

A first-principles approach is needed to justify why the principal components might be informative for simple, canonical noise models but also for identifying when middle components might be informative. This paper is a step in that direction In what follows, we provide a complete picture of how the spectrum of $X$ governs the informativeness of various components of the SVD of a data matrix $\wtX$ modeled as in (\ref{eq:sigplusnoise mat}).  To summarize our findings:
\begin{itemize}
\item The informative components correspond to isolated singular values that separate from the noise (or continuous looking) component of the spectrum,
\item Principal components are the most informative components when the noise (or the continuous looking) component of the spectrum is supported on one interval,
\item Middle components may be informative when the noise component of the spectrum is supported on multiple intervals,
\item Heterogeneities in the noise-only matrix can produce a disconnected noise spectrum,
\item It is possible for both principal and middle components to be informative and,
\item It is possible for the middle component to be informative even when the principal component is uninformative.
\end{itemize}

\begin{figure}[h]
\centering
\subfigure[Sample 1.]{
\includegraphics[trim = 100 250 100 250,clip = true, width = 1.75in]{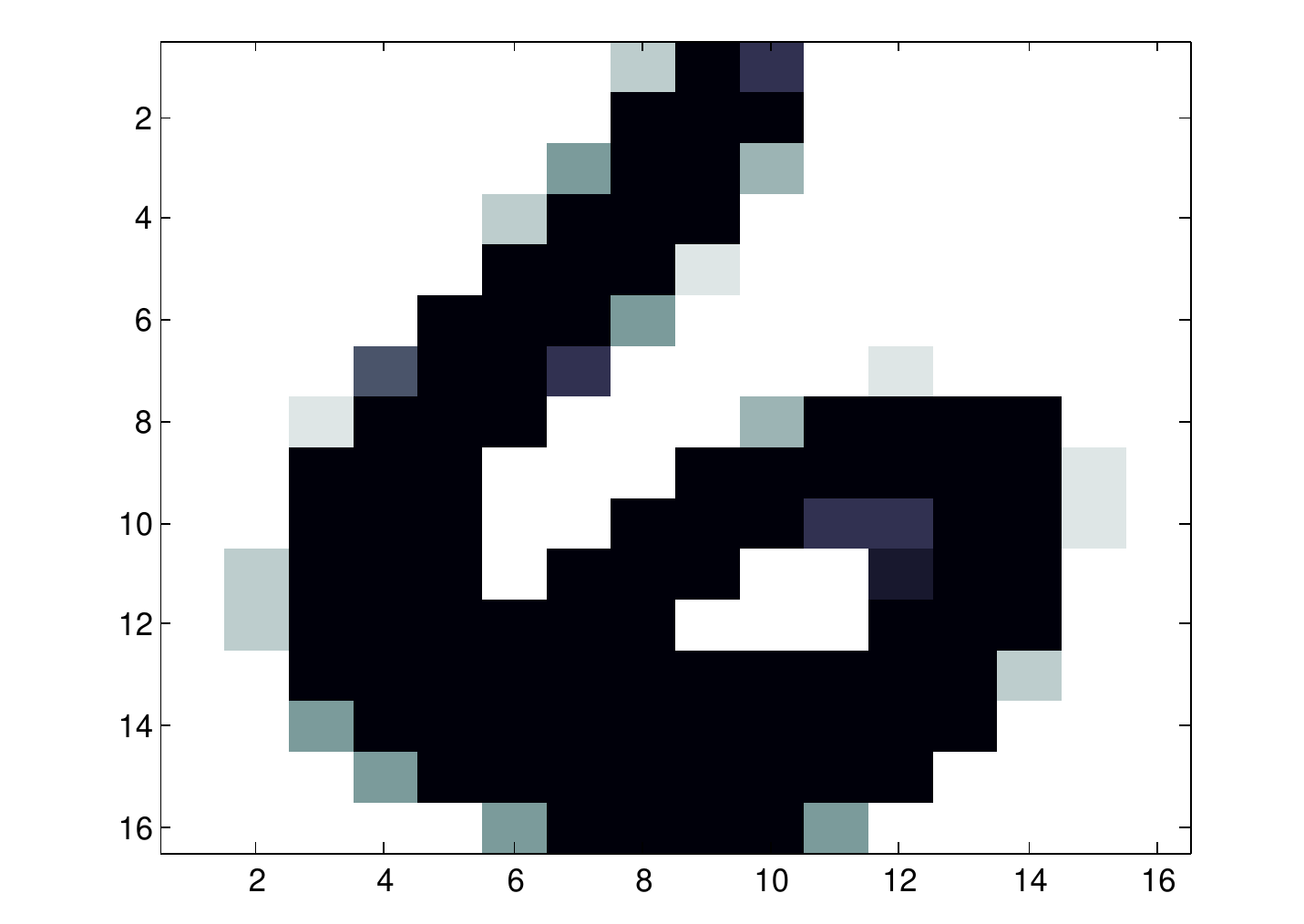}
}
\subfigure[Sample 2.]{
\includegraphics[trim = 100 250 100 250,clip = true, width = 1.75in]{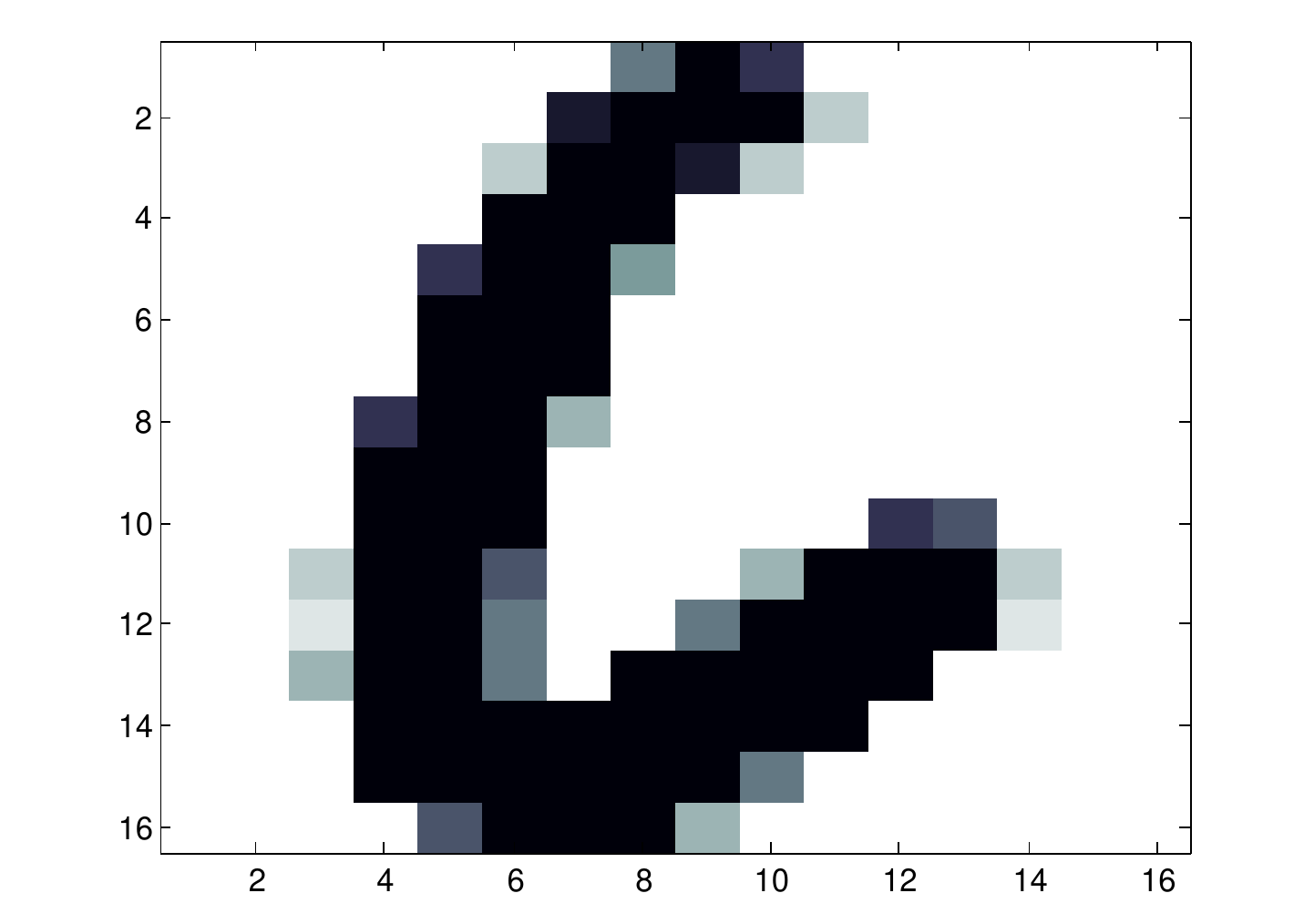}
}
\subfigure[Sample 3.]{
\includegraphics[trim = 100 250 100 250,clip = true, width = 1.75in]{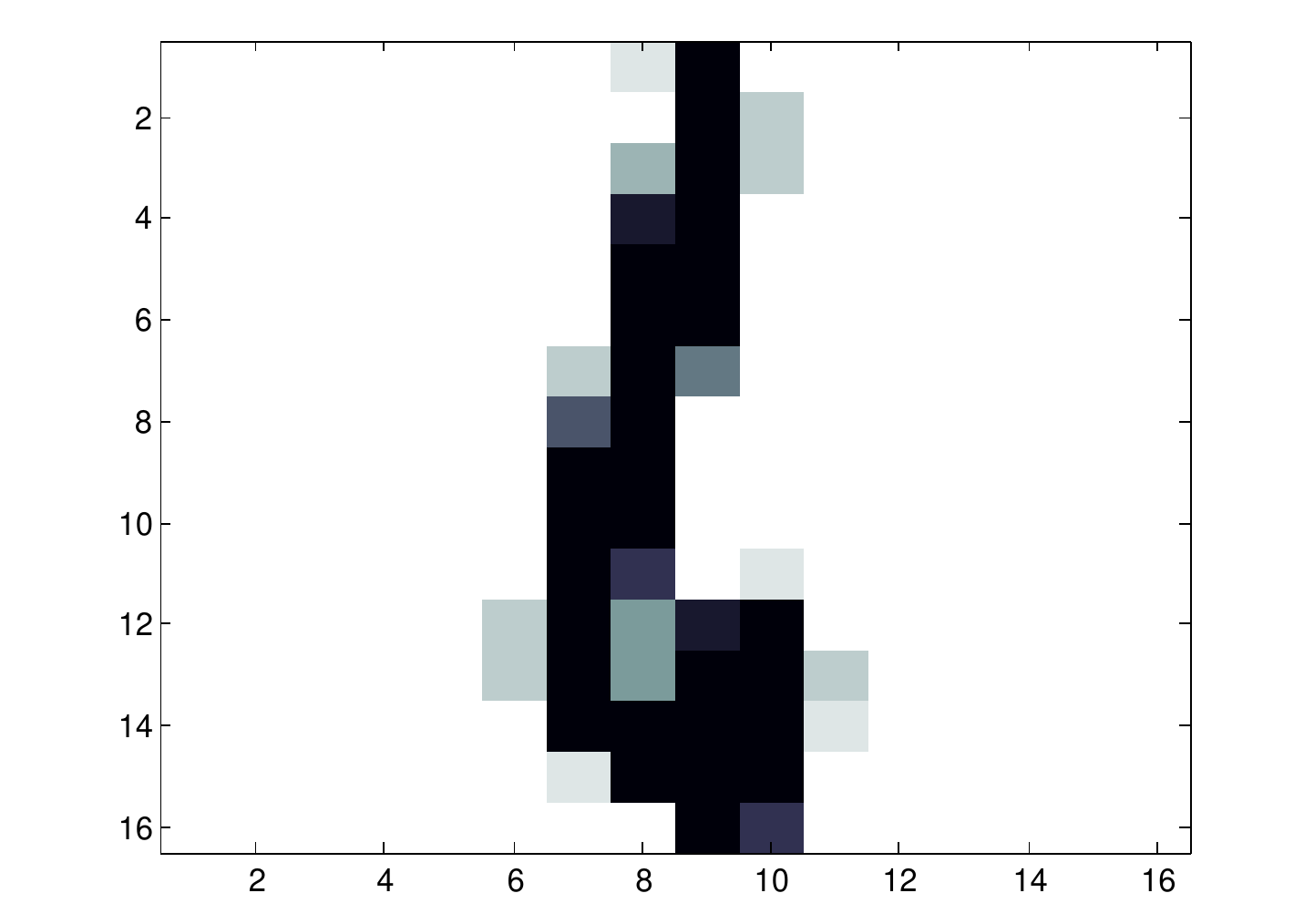}
}\\
\subfigure[The singular value spectrum of the training data matrix for the digit ``6'' is shown.]{
\centering
\includegraphics[scale = 0.85, trim = 0 170 0 170, clip = true, width=5.8in]{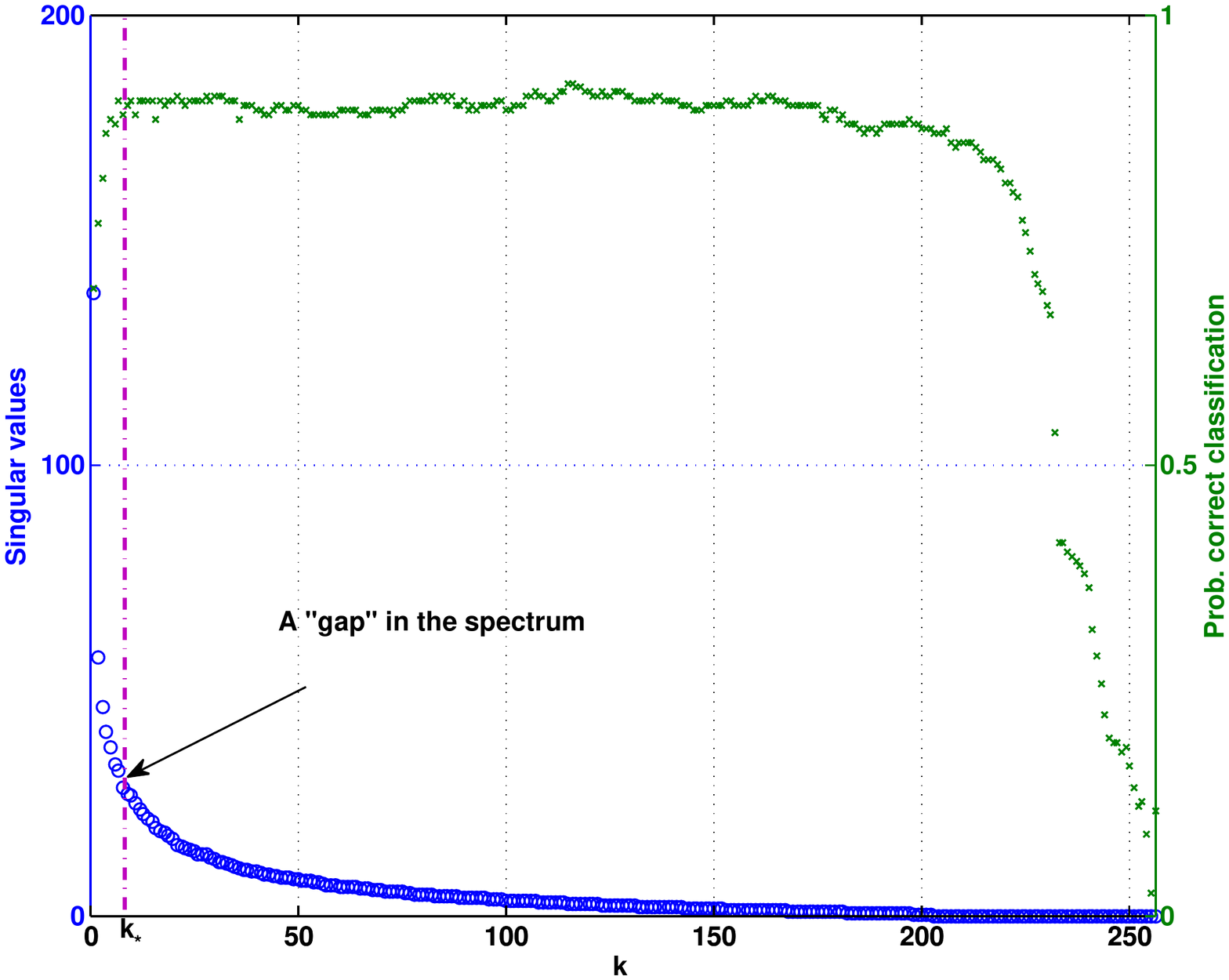}
}
\caption{(a) - (c) Three samples representing the digit ``6'' from the USPS handwritten digits database. Each $s$-pixel-by-$s$-pixel training image is converted into a  $n = s^2 \times 1$ column vector whose elements represent grayscale values. The training data matrix is formed by stacking the column vectors corresponding to every image in the labeled training data set alongside each other. (d)  displays on the singular values of the data matrix on the left axis. As in Figure \ref{fig:pca sv}, the singular value spectrum exhibits a continuous-looking portion (that may be interpreted as ``noise'') and a separated portion  (that may be interpreted as low-rank ``signal''). (right axis) A plot of a probability of correct classification versus $r$ plot where $r$ is the number of left singular vectors of the training set used for classification. \textit{Note that choosing $r$ based on the singular value ``gap separation'' heuristic yields near-optimal performance. }}
\label{fig:hand ex}
\label{fig:hand svd class}
\end{figure}

Our findings will allow the practitioner to better justify, by employing reasoning based on the \textit{entire} spectrum of $X$, when the use of principal components is warranted (as it is for the example in Figure \ref{fig:hand ex} ) and when the middle components might be more informative as in Figure \ref{fig:typical mca}. The next step in this line of inquiry, that is beyond the scope of this paper,  is the development of efficient computational methods for large data sets that can detect and extract informative middle components.

\begin{figure}[t]
\centering
\includegraphics[width=5.2in,trim = 65 220 65 220, clip   = true ]{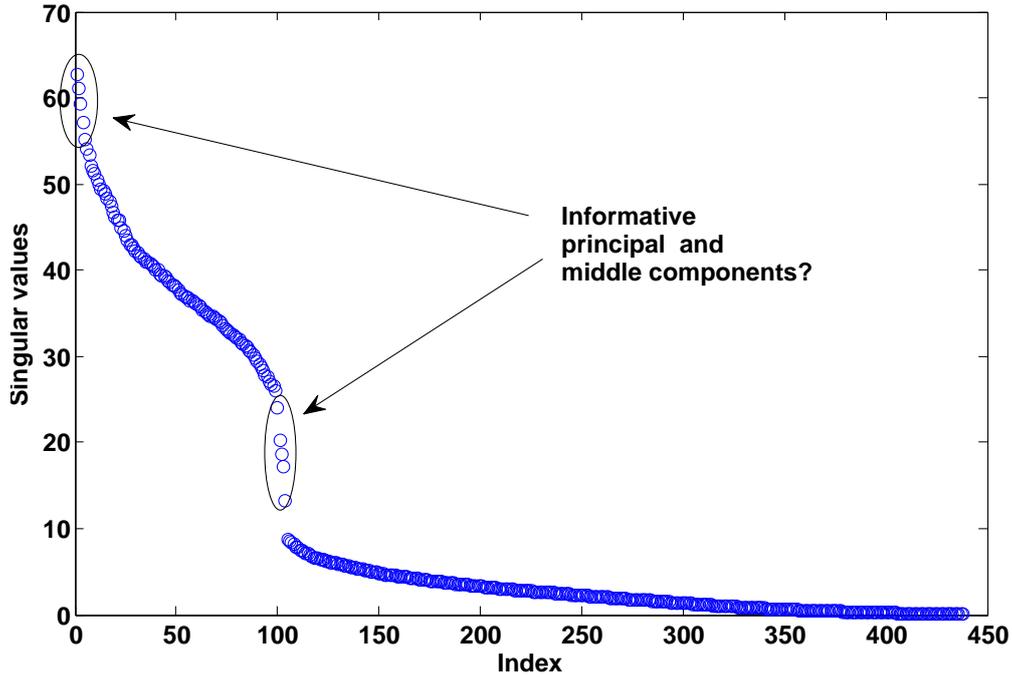}
\caption{A real-world data set with possibly informative principal and middle components.}
\label{fig:mca real world}
\end{figure}

We conclude by submitting Figure \ref{fig:mca real world} as evidence that our findings describe phenomena that might already be present in real-world data sets \footnote{ We thank Dr. James Preisig of the Woods Hole Oceanographic Institution for this dataset.} that might previously have been interpreted differently. Here we have a $438 \times 1200$ data matrix whose columns contains measurements made at a receiver sensor array and some of the past transmitted data symbols. The measurements were made over a time period where there were significant fluctuations in the noise levels.  The fluctuations in the channel transfer function constitute the low-rank ``signal'' here.

The plot of the singular values in  Figure \ref{fig:mca real world}  contains clusters of principal and middle eigenvalues that separate from the continuous looking portion of the spectrum. Our findings suggest that these are  informative principal and middle components. We hope that this work contributes to an increased understanding of the role played by the noise eigen-spectrum in shaping the informativeness of various SVD components and a recognition that there is much left to understand in terms of low-rank signal extraction from noisy data matrices.

We begin our exposition in Section \ref{sec:master equations} by examining how the spectrum of $\wtX$ is related to the spectrum of $X$.  We utilize the findings in Section \ref{sec:principal inform} to analyze a setting where the principal components are informative. In Section \ref{sec:middle inform} we describe a scenario when middle components can be informative while Section \ref{sec:main results} contains the main results which formalize the arguments presented in Sections \ref{sec:principal inform} and \ref{sec:middle inform}. We conclude in Section \ref{sec:noise models} with a discussion of which noise models can produce informative middle and principal components.

\clearpage

\section{The eigenvalues and eigenvectors of $\wtX$} \label{sec:master equations}
For expositional simplicity, let us consider the model in (\ref{eq:sigplusnoise mat}) with $r = 1$ and symmetric $\wtX$, so that
$$\wtX = S + X,$$
where $S = \theta \,uu^{*}$, for some arbitrary, non-random, unit norm column vector $u$.  We begin our investigation by examining how the eigenvalues and eigenvectors of $\wtX$ are related to the eigenvalues and eigenvectors of the low-rank signal matrix $S$.

Let $X = Q \Lambda Q^{*}$ be the eigen-decomposition of the noise-only random matrix $X$ (we have suppressed the subscript in $X_n$ for notational brevity), where $\Lambda=\diag(\la_1, \ldots, \la_n)$ and $Q$ are the eigenvectors of $X$. We assume that the noise-only random matrix $X$ is invariant, in distribution, under  orthogonal (or unitary) conjugation. This implies that the eigenvectors of $X$ are Haar-distributed and independent of its eigenvalues \cite[Th. 4.3.5]{hp00}. We will utilize this fact shortly.

\subsection{Eigenvalues of $\wtX$}\label{sec:max eig sketch}

The eigenvalues of $X+S$ are the solutions of the equation
$$ \det(zI - (X+S)) = 0 .$$
Equivalently, for $z $ such that $zI-X$ is invertible, we have $$zI-(X+S)=(zI-X) \cdot (I-(zI-X)^{-1}S),$$
so that
$$ \det(zI - (X+S)) = \det(zI - X) \cdot \det (I - (zI -X)^{-1}S).$$
Consequently, a simple argument reveals that the $z$ is an eigenvalue of $X+S$ and not an eigenvalue of $X$ if and only if $1$ is an eigenvalue of the matrix $(zI-X)^{-1}S$. But $(zI-X)^{-1}S= (zI-X)^{-1}\theta \,uu^{*}$ has rank one, so its only non-zero eigenvalue will equal its trace, which in turn is equal to $\theta u^{*}(zI-X)^{-1}u = \theta u^{*}Q(zI-\Lambda)^{-1}Q^{*}u $.

Let $v =Q^{*}u$. Then,  $z$ is an eigenvalues of $\wtX$ and not an eigenvalue of $X$ if and only if
\begin{equation}\label{eq:first master eqn}
\sum_{i=1}^{n} \dfrac{|v_{i}|^{2}}{z-\lambda_{i}} = \ff{\theta}.
\end{equation}
Let ${\mu_n}$ be the ``weighted" spectral measure of $X$, defined by \begin{equation}\label{eq:master eq mur1}\qquad
\mu_n = \sum_{i=1}^{n} | v_{i} |^{2} \delta_{\la_i} \qquad \textrm{ (the $v_i$'s are the coordinates of $v = Q^{*}u$).}
\end{equation}
Then any $z$ outside the  spectrum of $X$  is an eigenvalue of $\wtX$ if and only if
\begin{equation}\label{eq:master eq r1}
 \sum_{i=1}^{n}\dfrac{| v_i |^{2}}{z-\lambda_{i}}=:G_{\mu_n}(z)  = \ff{\theta},
\end{equation}
where $G_{\mu}(z)$ is the Cauchy transform of $\mu$ defined as
\begin{equation}\label{eq:cauchy transform}
G_{\mu}(z) = \int \dfrac{1}{z-x} d\mu(x).
\ee

Equation (\ref{eq:master eq r1}) describes the \textit{exact} relationship between the eigenvalues of $\wtX$ and the eigenvalues of $X$ and the dependence on the coordinates of the vector $v$ (via the measure $\mu_n$), which we will use shortly.

\subsection{Eigenvectors of $\wtX$}
Let $\widetilde{u}$ be a unit eigenvector of $X+S$ associated with the eigenvalue $ z $ that satisfies (\ref{eq:master eq r1}).  From the relationship $(X+S)\widetilde{u}=z \widetilde{u}$, we deduce that, for $S=\theta \,uu^{*}$,
$$\qquad\qquad
(zI-X)\widetilde{u}=S\widetilde{u}
=\theta uu^{*}\widetilde{u}
= (\theta u^{*}\widetilde{u}).u\qquad\qquad\textrm{(because $u^{*}\widetilde{u}$ is a scalar)},$$ implying that $\widetilde{u}$ is proportional to $(zI-X)^{-1}u$.

Since $\widetilde{u}$ has unit-norm, \be\label{250909.17h12}\widetilde{u}=\f{(zI-X)^{-1}u}{\sqrt{u^{*}(zI-X)^{-2}u}}\ee and \be\label{250909.17h05}  \uipvsq{u}{\wt{u}}=|u^{*}\widetilde{u}|^2=\f{(u^{*}Q(zI-\Lambda)^{-1}Q)^2}{u^{*}Q(zI-\Lambda)^{-2}Q^{*}u}=\f{G_{\mu_n}(z)^2}{\int  \f{\ud \mu_n(x)}{(z-x)^2}}=\ff{\theta^2\int  \f{\ud \mu_n(x)}{(z-x)^2}} .\ee
Notice that \be \int \f{\ud \mu_n(x)}{(z-x)^2} = - G'_{\mu_n}(z) \ee so that we have
\be \uipvsq{u}{\wt{u}} = - \ff{\theta^{2} G'_{\mu_n}(z)} \ee

Equation \eqref{250909.17h12} describes the relationship between the eigenvectors of $\widetilde{X}$ and the eigenvalues of $X$ and the dependence on the coordinates of the vector $v$ (via the measure $\mu_n$), which we will return to shortly.

\section{When principal components are the most informative components}\label{sec:principal inform}
We begin our investigation by considering a setting where the informative components do indeed correspond to the principal components. The picture we have developed so far is that the eigenvalues $z_{i}$ and the associated eigenvectors $\wt{u}_{i}$ of the signal-plus-noise data matrix $\wtX$ modeled as $\wtX = X + \theta uu^{*}$ satisfy the equations
\begin{subequations}\label{eq:eig picture}
\be
G_{\mu_n}(z_i) = \dfrac{1}{\theta},
\ee
\be
\label{eq:eigenvector inf}
|\langle \widetilde{u}_{i}, u \rangle|^{2} =  -\dfrac{1}{\theta^2}\cdot \dfrac{1}{G'_{\mu_n}(z_i)},
\ee
\end{subequations}
where
$$G_{\mu_n}(z) = \sum_{i=1}^{n}\dfrac{| v_i |^{2}}{z-\lambda_{i}}.$$

The expressions in (\ref{eq:eig picture}) provide insight on how the eigenvalues of $\wtX$ are related to the eigenvalues of $X$.

Figure \ref{fig:pca value} considers the $n = 5$ setting and shows how the expressions in (\ref{eq:eig picture}) provide insight on the informativeness of the eigenvalues and eigenvectors of $\wtX$.

By (\ref{eq:eig picture}a), the eigenvalues of $\wtX$ correspond to the values of $z$ where the horizontal line $1/\theta$ in Figure \ref{fig:pca value} intersects the curve $G_{\mu_n}(z)$.  Since   $G_{\mu_n}(z)$ has poles at the eigenvalues of $X$, all but the largest eigenvalue of $\wtX$ interlace the eigenvalues of $X$.  Consequently, $\lambda_{5} \leq \wt{\lambda}_{5} \leq \lambda_{4}$ and so on; there is no eigenvalue to the right of $\lambda_{1}$ and hence $\wt{\lambda}_{1}$ can be displaced by a greater amount, subject to $\wt{\lambda}_{1} - \lambda_{1} \leq \theta$.

Equation (\ref{eq:eig picture}b) reveals that the informativeness of an eigenvector, denoted by $\inform{i} : = \uipvsq{\wt{u}_{i}}{u}$, is inversely proportional to the negative slope of the function $G_{\mu_n}(z)$ evaluated at the eigenvalue $z_i = \wt{\lambda}_{i}$ of $\wtX$ associated with the eigenvector $\wt{u}_{i}$.

\subsection{Asymptotic analysis: Eigenvalues}
We now place ourselves in the high-dimensional setting. Let us assume that as $n \longrightarrow \infty$,
$$\mu_{X_n} = \dfrac{1}{n} \sum_{i=1}^{n} \delta_{\lambda_{i}}   \convas \mu_{X}, $$

where $\mu_X$ is a non-random probability and $\convas$ denotes almost sure convergence\footnote{The argument holds for other modes of convergence as well so we shall not explicitly specify the mode of convergence in the expository sections that follow.}. Assume that the largest and smallest eigenvalues of $X_n$ converge to $b$ and $a$, respectively and that $d \mu_X(z) > 0$ for all $z \in (a,b)$ so that the measure is supported on one connected interval. When  $X$ is a sample covariance matrix formed from a matrix with i.i.d.  Gaussian variables, the eigenvalues will satisfy this condition  \cite{silverstein1995empirical}.

The assumed convergence of the eigenvalues to a smooth limiting measure implies that as $n \to \infty$, if there were no signal, the eigenvalues would have a continuous looking spectrum as the spacing between successive eigenvalues goes to zero. By the same reasoning, when there is a signal, the picture developed in Figure \ref{fig:pca value}  says that all but the leading eigenvalue of $\wtX$ will be displaced insignificantly. Thus the $n-1$  eigenvalues will retain their continuous looking nature and will be tightly packed together.

As $n \to \infty$, only the largest eigenvalue will exhibit a significant $O(1)$ deviation relative to the corresponding eigenvalue in the noise-only setting (\textit{i.e.}, when $S = 0$). Since the second largest eigenvalue is also displaced insignificantly by a vanishing (with $n$) amount, this manifests as an $O(1)$ gap in the spectrum as in Figure \ref{fig:pca sv} and the use of the (principal) gap heuristic in (\ref{eq:khat}) for signal detection is justified.

We now investigate the fundamental limit of gap heuristic based signal detection. We first note that the vector $v = Q^{*}u$ is uniformly distributed on the unit hypersphere, and so, in the high-dimensional setting, $|v_{i}|^{2} \approx 1/n$ (with high probability) so that
$$ \sum_{i=1}^{n} |v_{i}|^{2} \delta_{\lambda_{i}} =: \mu_n  \approx \mu_{X} := \lim_{n \to \infty} \dfrac{1}{n} \sum_{i=1}^{n} \delta_{\lambda_{i}}. $$
A consequence of  $\mu_n  \to \mu_{X}$ is that $G_{\mu_n}(z) \to G_{\mu_{X}}(z)$. Inverting equation (\ref{eq:master eq r1}) after substituting these approximations yields the location of the largest  eigenvalue, in the $n \to \infty$ limit to be $G_{\mu_{X}}^{-1}(1/\theta)$.

Recall that we had assumed that the limiting probability measure of the noise-only random matrix $\mu_X$  is compactly supported on a single, connected interval $[a,b]$. Consequently, the Cauchy transform $G_{\mu_{X}}$ given by (\ref{eq:cauchy transform}) is well-defined for $z$ \textit{outside} $[a,b]$ and can tend to a limit $G_{\mu_X}(b^+)$ which may be bounded, \textit{i.e.} have $G_{\mu_X}(b^+) < +\infty$.

So long as $1/\theta  < G_{\mu_X}(b^+)$,  as in Figure \ref{fig:pca phase explanation}-(a),   we obtain $\lambda_{1}(\widetilde{X}) \approx G_{\mu_X}^{-1}(1/\theta) > b$. This results in an $O(1)$ gap between the largest eigenvalue and the edge of the spectrum and the gap heuristic will work. However, when $1/\theta \geq G_{\mu_X}(b^+)$, as in Figure \ref{fig:pca phase explanation}-(b), $\lambda_{1}(\widetilde{X}) \to  \lambda_{1}(X) = b$ and the gap heuristic will fail. To summarize:
\begin{center}
Principal gap based signal detection will asymptotically succeed iff $\theta > 1/G_{\mu_X}(b^+)$.
\end{center}

\subsection{Asymptotic analysis: Eigenvectors}

Recall our argument that since $Q$ is isotropically random, the vector $v = Q^{*}u$ is uniformly distributed on the unit hypersphere and $|v_{i}|^{2} \approx 1/n$ (with high probability) in the high-dimensional setting.  Consequently, we have  that $$G_{\mu_n}'(z) = - \dfrac{|v_i|^{2}}{(z-\lambda_i)^{2}} \approx -\dfrac{1}{n} \sum_{i=1}^{n} \dfrac{1}{(z-\lambda_{i})^{2}}.$$
Since all the $n$ eigenvalues of the noise-only matrix are concentrated on the connected interval $[a,b]$, the average spacing between the eigenvalues of $X$ is $O(1/n)$. Since the eigenvalues of $\wtX$ interlace the eigenvalues of $X$, $\wt{\lambda}_{i} - \lambda_{i} = O(1/n)$ for all but the largest eigenvalue. Hence $G_{\mu_n}'(z) = O(n)$ so that $\{\inform{i}\}_{i=2}^{n} =- 1/G'_{\mu_{n}}(z) |_{z = \widetilde{\lambda}_{i}} = O(1/n)$.

However,  $\wt{\lambda}_{1} - \lambda_{1} = O(1)$, so that $G_{\mu_n}'(\wt{\lambda}_{1}) = O(1)$ and $\inform{1} = O(1)$ implying that  the principal eigenvector is maximally informative  with a non-vanishing (with $n$) informativeness and the use of (\ref{eq:Shat}) in the estimation of $S$ is justified. We now investigate the fundamental limit of principal eigenvector based signal estimation.

In the asymptotic setting when $\mu_n \to  \mu_{X}$ and $z = \wt{\lambda}_{1} \to \rho$ we have that
$$\int  \f{\ud \mu_n(t)}{(z-t)^2}\to \int  \f{\ud \mu_X(t)}{(\rho-t)^2}=-G_{\mu_X}'(\rho),$$
so that
when $1/\theta <G_{\mu_X}(b^+)$, which implies that $\rho>b$,    we have
 $$|\lan \widetilde{u}_{1}, u\ran|^2\convas
 \ff{\theta^2\int  \f{\ud \mu_X(t)}{(\rho-t)^2}}=\f{-1}{\theta^2G_{\mu_X}'(\rho)}>0,$$
 whereas when $1/\theta \ge G_{\mu_X}(b^+)$ and if $\mu_X$ is such that  $G_{\mu_X}$ has infinite derivative at $\rho=b$, we have
$$|\lan \wt{u}_1, u \ran | \convas 0. $$
Hence when $\theta \leq 1/G_{\mu_X}(b^+)$ and if $G'_{\mu_X}(b^+) = \infty$, then the all components have vanishing (with $n$) informativeness.
To summarize, when $\theta > 1/G_{\mu_X}(b^+)$:
\begin{quote}
Principal components are the most informative components when the noise eigen-spectrum is contained on a single, connected interval .
\end{quote}
The eigen-spectrum of a Wishart distributed sample covariance matrix with identity covariance satisfies this condition. It is thus a happy coincidence that principal components are the most informative components for the simplest noise matrix model. We now consider the setting where the noise eigen-spectrum is (asymptotically) supported on multiple disconnected intervals.

\begin{figure}
\centering
\includegraphics[trim = 75 125 0 0, clip = true,width=5.55in]{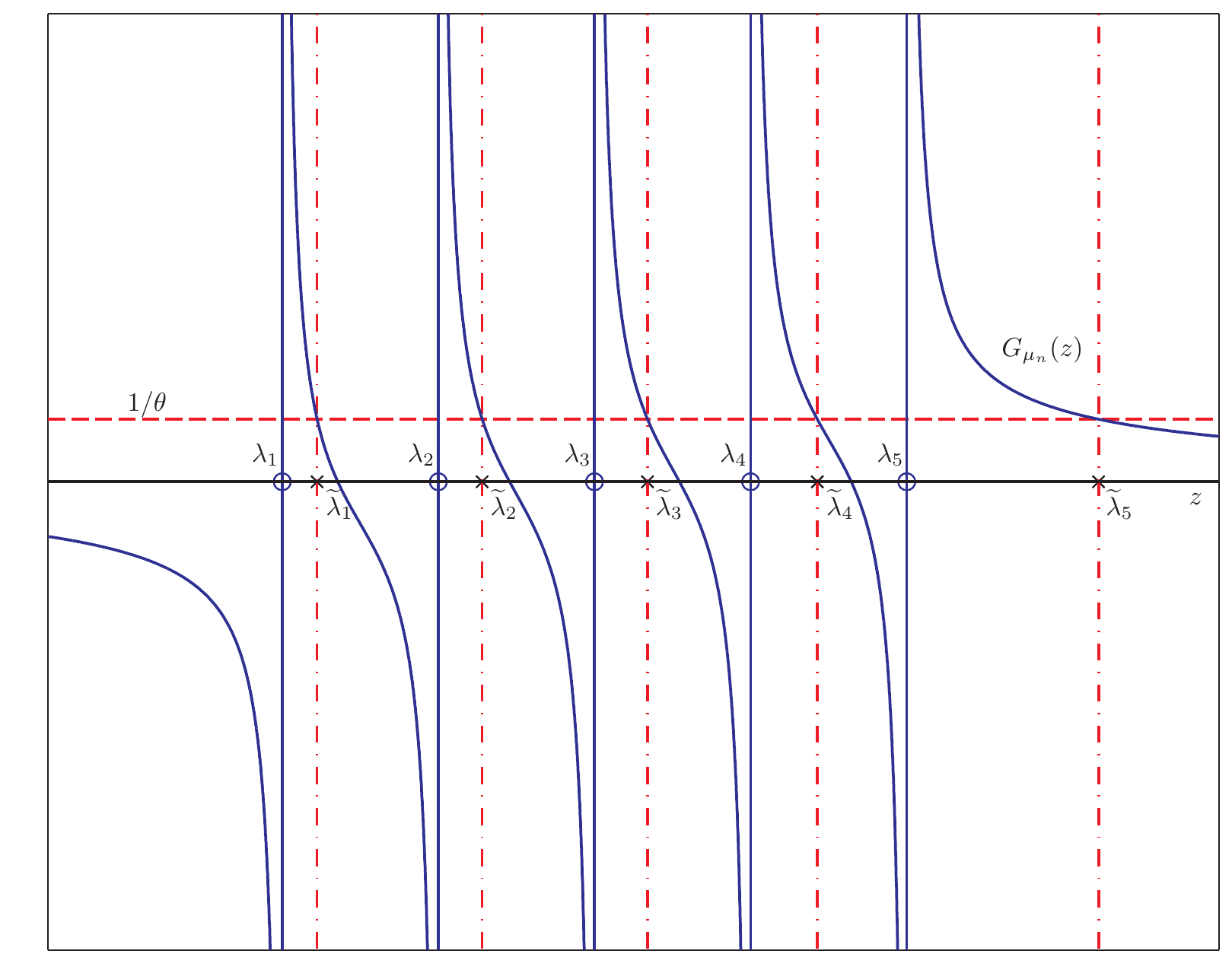}
\caption{The relationship in (\ref{eq:eig picture}a) between the eigenvalues of $\wtX= \theta uu^* + X$ and the eigenvalues of $X$ is depicted here. Notice the interlacing of the bulk eigenvalues and the emergence of the principal eigen-gap.}
\label{fig:pca value}
\end{figure}

\begin{figure}
\centering
\subfigure[When $1/\theta < G_{\mu}(b)$, $\lambda_{1}(\widetilde{X}) \to \rho = G^{-1}_{\mu}(\tfrac{1}{\theta})$ and there is a principal eigen-gap.]{
\includegraphics[trim =130 315  80 320, clip = true, width = 5.5in]{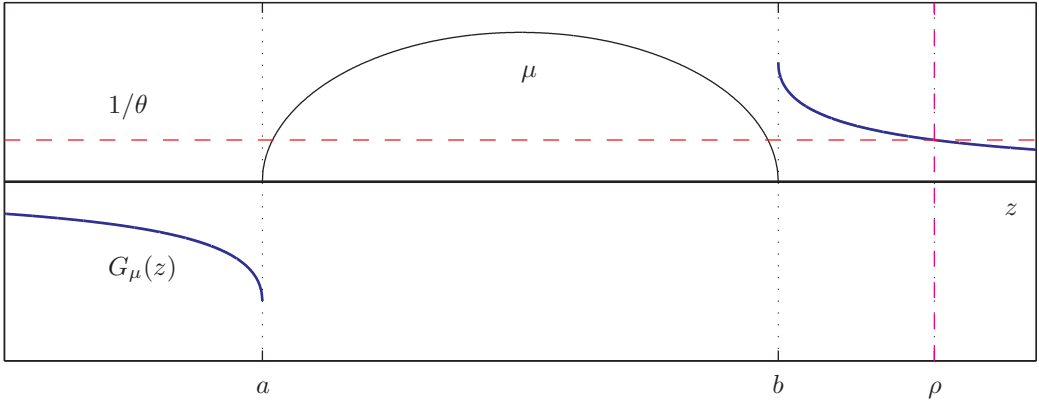}
}\\
\subfigure[When $1/\theta =  G_{\mu}(b)$, $\lambda_{1}(\widetilde{X}) \to b$ and there is no principal eigen-gap.]{
\psfrag{d}{$1/\theta$}
\includegraphics[trim =130 315  80 320, clip = true, width = 5.5in]{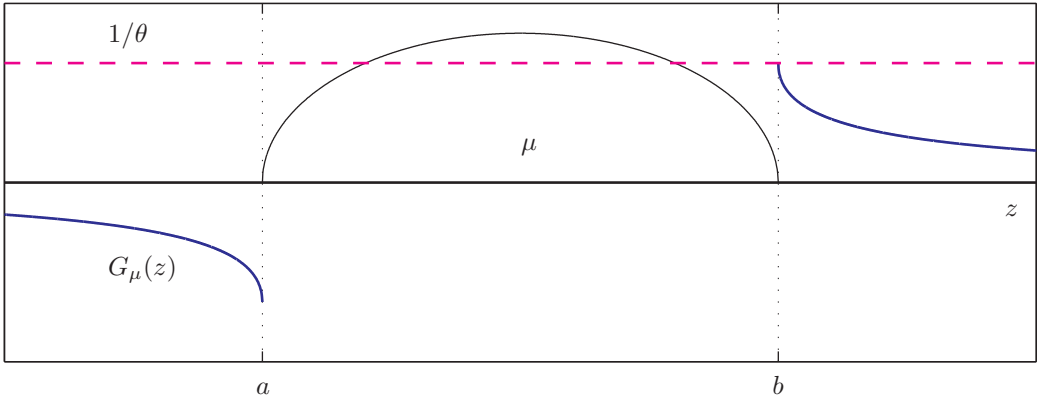}
}
\caption{The evolution of the informativeness of the principal eigen-gap for different values of $\theta$. In (a), where $\theta$ is large, the principal eigen-gap is informative; (b) the principal eigen-gap vanishes when $1/\theta = G_{\mu}(b^+)$ and the signal is undetectable using principal eigen-gap based methods.}
\label{fig:pca phase explanation}
\end{figure}

\clearpage
\section{When middle components are informative}\label{sec:middle inform}

\subsection{Asymptotic analysis: eigenvalues}
Consider a setting where the eigenvalues of the noise-only matrix $X$ are supported on multiple intervals as in Figure \ref{fig:typical mca}. This corresponds to letting the $\mu_X$ obtained in the $n \longrightarrow \infty$ limit of $\mu_{X_n}$  being supported on $\ell = 2$ intervals. Thus we may model $\mu_X$ as
$$\mu_X = p_1 \mu_{1} + p_2  \mu_{2},$$
where $p_1 + p_2 = 1$. Here $p:=p_1 \in (0,1)$ and the measures $\mu_{1}$ and $\mu_2$ are non-random probability measures supported on $[a_{1},b_{1}]$ and $[a_2,b_2]$, respectively with $\ud\mu_{i}(z) > 0$ for $z \in (a_{i},b_{i})$ for $i = 1, 2$.  We suppose that $a_{2} < b_{2} < a_{1} < b_{1}$, as depicted in Figure \ref{fig:mca phase explanation}-(a). For $k = 0, 1, 2$, define $c_{j} = \sum_{i=0}^{j} p_{i}$ with $p_0 := 0$, $c_0 := 0$ and $c_{2} := 1$. We assume that  for $j = 1, 2$, $\lambda_{n c_{j-1} + 1} \convas b_{j}$ and that $\lambda_{n c_{j}} \convas a_{j}$. For expositional simplicity we assume that $n c_j$ is an integer.  When $X$ is a sample covariance matrix formed from a matrix with Gaussian entries  having a covariance matrix with an adequately-separated covariance eigen-spectrum then the sample eigenvalues will satisfy this condition. Section \ref{sec:noise models} contains additional examples and elaborates on when the covariance eigen-spectrum in separated enough.

The assumed convergence of the eigenvalues to a smooth limiting measure implies that as $n \to \infty$, if there were no signal, the eigenvalues would have a continuous looking spectrum as the spacing between successive eigenvalues goes to zero.

By the same reasoning, when there is a signal, the picture developed in Figure \ref{fig:pca value} when adapted as in Figure \ref{fig:mca value} for the disjoint interval setting (here $\ell = 2$) reveals (via (\ref{eq:eig picture}) that the leading  eigenvalue $\widetilde{\lambda}_{1}$ \textbf{and} an additional, middle, eigenvalue $\widetilde{\lambda}_{n c_1+1}$ will exhibit a significant $O(1)$ deviation relative to the corresponding eigenvalue in the noise only setting. The middle eigenvalue emerges from the bulk spectrum because $\lambda_{nc_1} - \lambda_{nc_1 + 1} = O(1)$ as $n \to \infty$.

The remaining $n - 2$ eigenvalues will be displaced insignificantly and will remain tightly packed together,  thereby retaining their continuous looking appearance. Consequently, there will be two O(1) eigen-gaps in the spectrum betraying the presence of a low-rank signal. Thus, here too, the use of the gap heuristic is justified.

The emergence of an informative middle eigenvalue in this setting due to the presence of a large gap in the noise eigen-spectrum may be viewed as a form of aliasing.

We now investigate the fundamental limit of  gap heuristic based signal detection so we might understand when not accounting for the middle eigen-gap might lead to suboptimal detection performance.

As before, we note that the vector $v = Q^{*}u$ is uniformly distributed on the unit hypersphere, and so, in the high-dimensional setting, $|v_{i}|^{2} \approx 1/n$ (with high probability) so that
$$ \sum_{i=1}^{n} |v_{i}|^{2} \delta_{\lambda_{i}} =: \mu_n  \approx \mu_{X} := \lim_{n \to \infty} \dfrac{1}{n} \sum_{i=1}^{n} \delta_{\lambda_{i}}. $$
A consequence of  $\mu_n  \to \mu_{X}$ is that $G_{\mu_n}(z) \to G_{\mu_{X}}(z)$. Inverting equation (\ref{eq:master eq r1}) after substituting these approximations yields the location of the largest  eigenvalue, in the $n \to \infty$ limit to be $G_{\mu_{X}}^{-1}(1/\theta)$. This results in multiple (i.e. principal and middle) eigenvalues that separate from the bulk spectrum precisely when the functional inverse is multi-valued (for a domain outside the region of support).

Recall our assumption that the limiting probability measure of the noise-only random matrix $\mu_X$ x is compactly supported on $\ell=2$ disjoint intervals $\{[a_i,b_i]\}_{i=1}^{2}$. Consequently, the Cauchy transform $G_{\mu_{X}}$ given by (\ref{eq:cauchy transform}) is well-defined for $z$ \textit{outside} $\cup_{i=1}^{2} [a_i,b_i]$ and is strictly decreasing with increasing $z$ on open intervals  $\left(\cup_{i=1}^{2} [a_i,b_i] \right)^{c} $ outside the support of $\mu_X$, as depicted in Figure \ref{fig:mca phase explanation}.

Thus so long as $1/\theta < G_{\mu_X}(b_1^+)$, $\lambda_{1} \to \rho_1>b_1$ and an $O(1)$ principal eigen-gap will manifest. Conversely, if $1/\theta \geq  G_{\mu_X}(b_1^+)$, as in Figure \ref{fig:mca phase explanation}-(b), $\lambda_{1} \to b_1$ and there will be no principal eigen-gap.

Similarly,  if $1/\theta < G_{\mu_X}(b_2^+)$ and $1/\theta > G_{\mu_X}(a_1^-)$, $\lambda_{nc_1 +1 } = \lambda_{n p + 1} \to \rho_2>b_2$ and an $O(1)$ middle eigen-gap will manifest. Conversely, , as shown in Figure \ref{fig:mca phase explanation}-(c),  if $1/\theta > G_{\mu_X}(b_2^+)$ then $\lambda_{nc_1 +1 } = \lambda_{n p + 1} \to b_2$ and there will be no $O(1)$ middle eigen-gap. However when $1/\theta < G_{\mu_X}(a_1^-)$, then $\lambda_{nc_1 +1 } = \lambda_{n p + 1} \to a_1$  and technically speaking there is an $O(1)$ middle eigen-gap except that this gap is indistinguishable from the gap in the spectrum that appears even when there is no signal.
 
Thus principal eigen-gap based signal detection for weak signals (or small $\theta$) fails whenever $\theta < 1/ G_{\mu_X}(b_{1}^+)$ while middle eigen-gap detection fails whenever  $\theta < 1/ G_{\mu_X}(b_{2}^+)$. If $G_{\mu_X}(b_2^+) > G_{\mu_X}(b_1^+)$, as depicted in Figure \ref{fig:mca phase explanation}, then a weak signal that is undetectable using the principal eigen-gap heuristic would have remained detectable if the middle eigen-gap were considered. This is why the middle eigen-gap in Figure \ref{fig:typical mca} was informative while the principal eigen-gap was not.  In such settings, detection using \textbf{only} the principal eigen-gap detection is suboptimal.

\subsection{Asymptotic analysis: eigenvectors}

Equation (\ref{eq:eigenvector inf}) reveals that  the informativeness of an eigenvector $\widetilde{u}_{i}$, relative to the signal eigenvector $u$ is given by the expression
$$|\langle \widetilde{u}_{i}, u \rangle|^{2} =  -\dfrac{1}{\theta^2}\cdot \dfrac{1}{G'_{\mu_n}(z_i)},$$
where
$$G_{\mu_n}'(z) = - \dfrac{|v_i|^{2}}{(z-\lambda_i)^{2}} \approx -\dfrac{1}{n} \sum_{i=1}^{n} \dfrac{1}{(z-\lambda_{i})^{2}}.$$
The eigenvalues of the noise-only matrix are concentrated on the disjoint intervals $[a_1,b_1]$ and $[a_2,b_2]$. Thus, the average spacing between the successive eigenvalues of $X$ within each interval is $O(1/n)$. Since the eigenvalues of $\wtX$ interlace the eigenvalues of $X$, $\wt{\lambda}_{i} - \lambda_{i} = O(1/n)$ for all but the largest eigenvalue and the middle eigenvalue as in Figure \ref{fig:mca phase explanation}-(a). Hence $G_{\mu_n}'(z) = O(n)$ so that $\{\inform{i}\}_{i=2}^{n} = 1/G'_{\mu_{n}}(z) |_{z = \widetilde{\lambda}_{i}} = O(1/n)$.

As before, we note that  $\wt{\lambda}_{1} - \lambda_{1} = O(1)$   so that $G_{\mu_n}'(\wt{\lambda}_{1}) = O(1)$ and $\inform{1} = O(1)$ implying that the principal eigenvector is informative  with a non-vanishing (with $n$) informativeness and the use of (\ref{eq:Shat}) in the estimation of $S$ is justified. However, what emerges from the picture in Figure \ref{fig:mca value} is that since  $\lambda_{np}- \lambda_{np+1} = O(1)$ we have that $\wt{\lambda}_{np+1} - \lambda_{np+1} = O(1)$  and by the same argument, $\inform{np+1} = O(1)$ as well. Thus the middle eigenvector associated with the middle eigenvalue that exhibits an eigen-gap is also informative. Employing it in the estimation of $S$ in (\ref{eq:Shat}) would improve estimation performance.

\begin{figure}[t]
\centering
\vspace{-0.1cm}
\includegraphics[trim = 130 190 50 220, clip = true,width=5.5in]{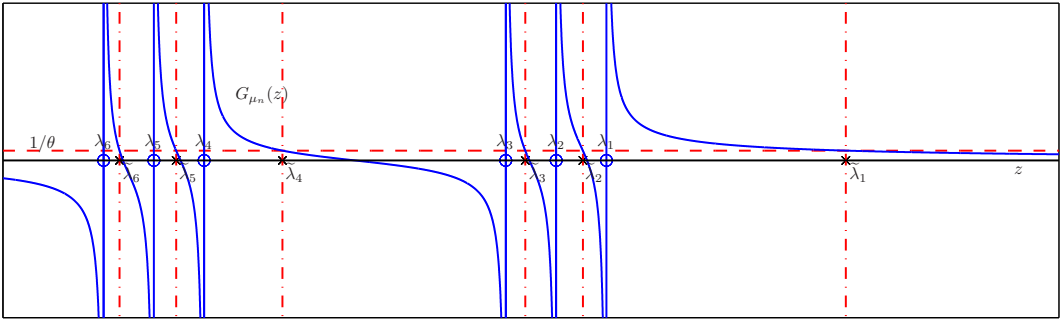}
\caption{The relationship in (\ref{eq:eig picture}a) between the eigenvalues of $\wtX= \theta uu^* + X$ and those of $X$ is depicted here when the eigenvalues of $X$ are supported on two $O(1)$ separated intervals. Notice the interlacing of the bulk eigenvalues and the emergence of the principal and middle eigen-gap . Contrast this to Figure \ref{fig:pca value} when only the principal eigen-gap emerges.}
\label{fig:mca value}
\vspace{-0.1cm}
\end{figure}

Extending this argument further, in the $n \to \infty$ limit, when $\mu_n \to  \mu_{X}$ suppose $\mu_X$ is such that $G'_{\mu_X}(z)=- \infty$ for $z = b_1, b_2$. Then we have that whenever $1/\theta < G_{\mu_X}(b_{1}^+)$,
 $$|\lan \widetilde{u}_{1}, u\ran|^2\convas
 \ff{\theta^2\int  \f{\ud \mu_X(t)}{(\rho_1-t)^2}}=\f{-1}{\theta^2G_{\mu_X}'(\rho_1)}>0,$$
but if $1/\theta \geq G_{\mu_X}(b_{1}^+)$ as in Figure \ref{fig:mca phase explanation}-(b),(c), then $|\lan \wt{u}_1, u \ran | \convas 0,$
and the principal component becomes uninformative. Employing the same argument for the middle eigenvector reveals that so long as $G_{\mu_X}(a_1^-) \leq 1/\theta < G_{\mu_X}(b_2^+)$ then $\widetilde{\lambda}_{np+1} \convas \rho_2$ and the corresponding eigenvector is informative \textit{i.e.},
 $$|\lan \widetilde{u}_{np+1}, u\ran|^2\convas \f{-1}{\theta^2G_{\mu_X}'(\rho_2)}>0.$$
When $1/\theta \geq G_{\mu_X}(b_{2}^+)$ as in Figure \ref{fig:mca phase explanation}-(c), then
$$|\lan \wt{u}_1, u \ran | \convas 0,$$
and the middle component becomes uninformative. Evidently, if $1/G_{\mu_X}(b_2^+) > 1/G_{\mu_X}(b_1^+)$ as in Figure \ref{fig:mca phase explanation} then the middle eigenvector will stay informative for a regime of small $\theta$ where the principal eigenvector is uninformative. More generally, if both the principal and the middle eigenvectors are informative then principal eigenvector will be more informative if $-1/G'_{\mu_X}(\rho_1) > -1/G'_{\mu}(\rho_2)$ and vice versa. This is determined by the structure of the noise spectrum. To summarize:
\vspace{-0.1cm}
\begin{itemize}
\item Principal gap based signal detection will asymptotically succeed iff $\theta > 1/G_{\mu_X}(b_1^+)$,
\item Middle gap based detection will asymptotically succeed despite principal gap based detection failing whenever $G_{\mu_X}(b_2^+) > G_{\mu_X}(b_1^+)$.
\item The eigenvectors associated with principal or middle eigenvalues that exhibit an eigen-gap will be informative
\item The eigenvectors will be uninformative when the eigen-gap vanishes.
\end{itemize}
The emergence of informative middle eigenvalues and eigenvector whenever there is a gap in the noise eigen-spectrum may be viewed as a form of signal (subspace) aliasing.

\begin{figure}[h]
\centering
\subfigure[When $1/\theta < G_{\mu_X}(b_1^+)$ and $1/\theta < G_{\mu_X}(b_2^+)$, $\widetilde{\lambda}_{1}(\widetilde{X}) \to \rho_1 > b_1$ and $\widetilde{\lambda}_{n p +1} \to \rho_2 > b_{2}$. ]{
\includegraphics[trim =130 310  80 320, clip = true, width = 6.25in]{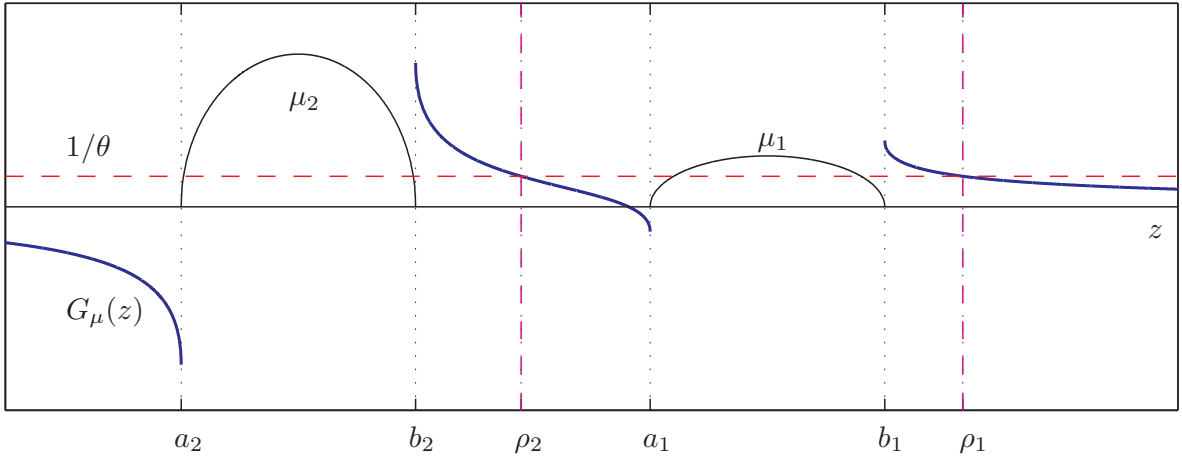}
}\\[-0.15cm]
\subfigure[When $1/\theta = G_{\mu_X}(b_{1}^+)$, $\wt{\lambda}_1 \to b_1$. Since $1/\theta < G_{\mu_X}(b_2^+)$, $\widetilde{\lambda}_{n p+1} \to \rho_2 >b_2$.]{
\includegraphics[trim =130 310  80 320, clip = true, width = 6.25in]{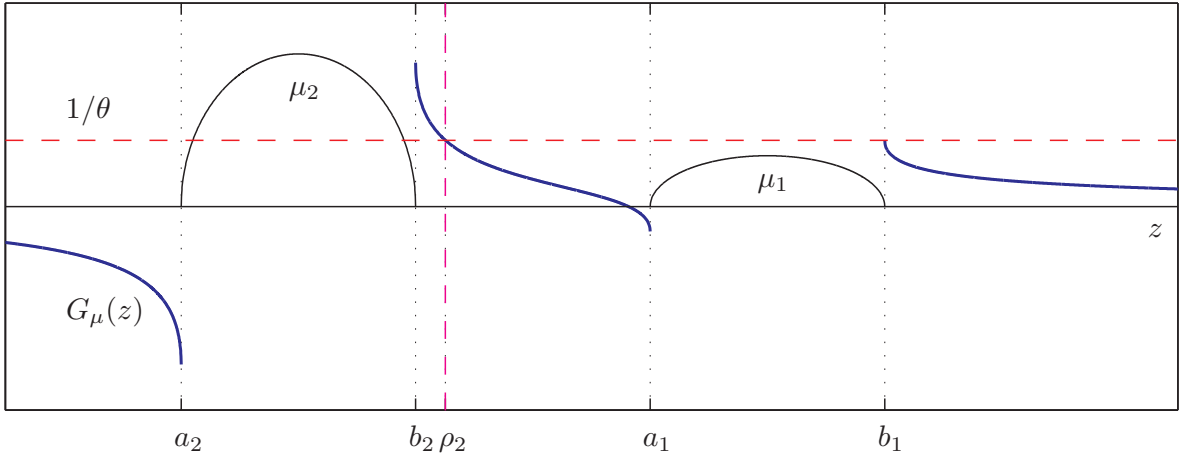}
}
\\[-0.15cm]
\subfigure[When $1/\theta >G_{\mu_X}(b_{1}^+)$, $\wt{\lambda}_1 \to b_1$. Since $1/\theta = G_{\mu_X}(b_2^+)$, $\wt{\lambda}_{n p+1} \to b_2$.]{
\psfrag{d}{$1/\theta$}
\includegraphics[trim =130 310  80 320, clip = true, width = 6.25in]{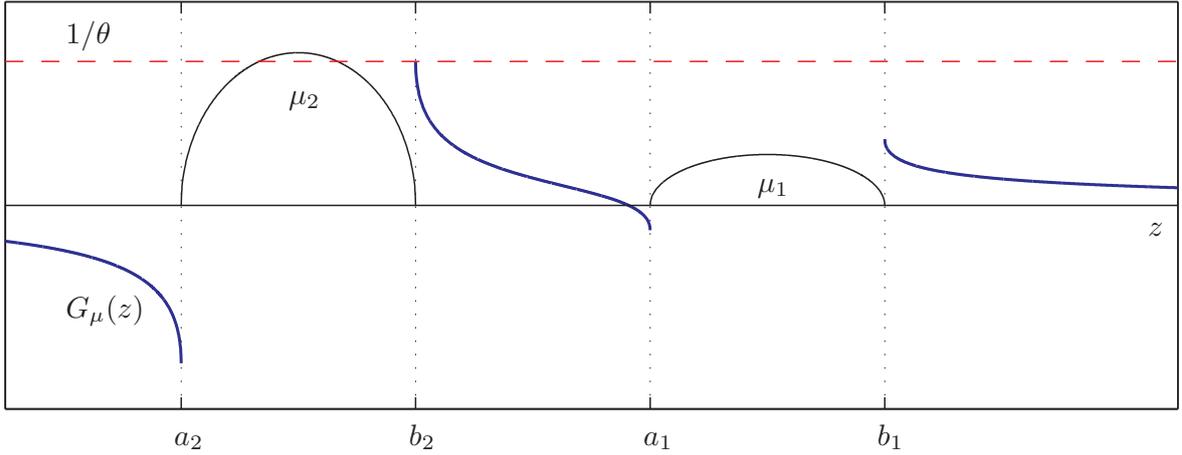}
}
\vspace{-0.35cm}
\caption{The evolution of the informativeness of the principal and the middle eigen-gaps for different values of $\theta$. In (a), where $\theta$ is large, both the principal and middle eigen-gaps are informative; (b) the principal eigen-gap vanishes but the middle eigen-gap persists; (c) both the principal and middle eigen-gaps vanishes and the signal is undetectable using eigen-gap based methods. The important point to note here is that the middle eigen-gap reveals the presence of a signal even when the principal eigen-gap does not.}
\label{fig:mca phase explanation}
\end{figure}

\clearpage

\section{Main results}\label{sec:main results}

\subsection{Eigenvalues and Eigenvectors}\label{sec:main ev}
Let  $X_n$ be  an $n \times n$ symmetric (or Hermitian) random matrix whose ordered eigenvalues we denote by $\la_{1}(X_n) \geq \cdots \geq \la_{n}(X_n)$. Let $\mu_{X_{n}}$ be the empirical eigenvalue distribution, \ie, the probability measure defined as
\[
\mu_{X_{n}} = \frac{1}{n} \sum_{i=1}^{n} \delta_{\la_{i}(X_{n})}.
\]
Assume that the probability measure $\mu_{X_{n}}$ converges almost surely weakly, as $n \longrightarrow \infty$, to a non-random compactly supported  \pro measure $\mu_X$ that is supported on $\ell$ disjoint intervals so that
$$\mu_X = \sum_{j=1}^{\ell} p_{j}\, \mu_{j} ,$$
where for $j = 1, \ldots, \ell$, the measures $\mu_{j}(x)$ is a non-random probability measure are supported on $[a_{j},b_{j}]$ with  $d\mu_{j}(z) > 0$ for $z \in (a_{j},b_{j})$ and $a_{\ell} < b_{\ell} < a_{\ell - 1} < b_{\ell-1} < \ldots  < a_{1} < b_{1}$. Define $c_{k} = \sum_{i=0}^{k} p_{i}$ with $p_0 := 0$, $c_0 := 0$ and $c_{\ell} := 1$. We assume that  for $j = 1, \ldots,  \ell$, $\lambda_{\lceil n c_{j-1} \rceil + 1} \convas b_{j}$ and that $\lambda_{\lceil n c_{j} \rceil} \convas a_{j}$, where $\lceil n c_{j-1} \rceil$ denotes the smallest integer greater than or equal to $n c_{j-1}$.

For a given $r \ge 1$, let $\tta_1\ge\cdots\ge \tta_r$ be deterministic non-zero real numbers, chosen independently of $n$.  For every $n$, let $P_n$ be an $n\times n$ symmetric (or Hermitian) random matrix   having rank $r$ with its $r$ non-zero eigenvalues equal to $\theta_{1}, \ldots, \theta_{r}$.

Recall that a symmetric (or Hermitian) random matrix is said to be {\it orthogonally invariant} (or {\it unitarily invariant}) if its distribution is invariant under the action of the orthogonal (or unitary) group under conjugation.

We suppose that $X_n$ and $P_n$ are independent and that $X_n$, the noise-only, matrix  is unitarily invariant while the low-rank signal matrix $P_n$ is non-random.

\subsubsection{Notation} Throughout this paper, for $f$ a function and $c\in \R$, we set $$f(c^+):=\lim_{z\downarrow c}f(z)\,;\qquad f(c^-):=\lim_{z\uparrow c}f(z),$$
we also let $\convas$ denote almost sure convergence. The ordered eigenvalues of an $n\times n$ Hermitian matrix $M$ will be denoted by $\la_1(M)\ge\cdots \ge \la_n(M)$. Lastly, for a subspace $F$ of a Euclidian space $E$ and a vector $x\in E$, we denote the norm of the orthogonal projection of $x$ onto $F$ by $\lan x, F\ran$.

Consider the rank $r$ additive perturbation of the random matrix $X_n$ given by
$$\wtX = X_n+P_n.$$
For this model, we establish the following results.

\begin{Th}[Eigen-gap phase transition]\label{140709.main}
The eigenvalues of $\wtX $ exhibit the following    behavior as $n \longrightarrow \infty$. We have that for each $1\le  i \leq r$ and $ 1 \leq j \leq \ell$,
$$\la_{\lceil n c_{j-1}\rceil  + i}(\wtX) \convas \begin{cases} G_{\mu_X,(b_j,a_{j-1})}^{-1}(1/\theta_i)& \textrm{ if } 1/G_{\mu_X}(b_{j}^+) < \theta_i < 1/G_{\mu_X}(a_{j-1}^{-}),\\ \\
b_{j} &\textrm{ if } \theta_i < 1/G_{\mu_X}(b_j^+),\\ \\
a_{j-1} & \textrm{ if } \theta > 1/G_{\mu_X}(a_{j-1}^-)
\end{cases}$$
Here,
$$G_{\mu_X}(z)=\int\f{1}{z-t}\ud\mu_X(t) \qquad \textrm{for } z \notin \supp \mu_X,$$
is the Cauchy transform of $\mu_X$,  $G_{\mu_X,(b_j,a_{j-1})}^{-1}(\cdot)$ is its functional inverse for $G_{\mu_X}(z)$ for $z \in (b_j,a_{j-1})$ and $a_{0} := + \infty$.
\end{Th}
\begin{proof}
The result is obtained by following the approach taken in \cite[pp. 511-514]{benaych2011eigenvalues} for proving Theorem 2.1. The key difference is that we are explicitly considering measures $\mu_X$ supported on multiple (disconnected) intervals so that the Cauchy transform of $\mu_X$ can have multiple inverses as in Figure \ref{fig:mca phase explanation}. For those values of $\theta$ such that $G_{\mu_X}^{-1}(1/\theta)$  is multi-valued, as many eigenvalues of $\wtX$ as there are values of $z$ such that $z = G_{\mu_X}^{-1}(1/\theta)$ will exhibit the eigen-gaps identified.
\end{proof}

\begin{Th}[Informativeness of the eigenvectors]\label{180709.13h39} Assume throughout that $\theta >0$ and let $G_{\mu_X}(a_0^-) = +\infty$. Consider   ${i_0}\in \{1, \ldots, r\}$  \st $1/\theta_{i_0} \in \bigcup_{j = 1}^{l} (G_{\mu_X}(a_{j-1}^-),G_{\mu_X}(b_{j}^+) )$. For each such $i_0$, consider $j(i_0) \equiv j \in \{1, \ldots, l\}$ such that $1/\theta_{i_0} \in (G_{\mu_X}(a_{j-1}^{-}),G_{\mu_X}(b_j^+) )$ and let $\widetilde{u}$ be a unit-norm eigenvector of $\wtX$  associated with the eigenvalue $\widetilde{\lambda}_{\lceil  nc_{j-1} \rceil +i_0} $. Then we have, as $n\lto\infty$,

\flushleft (a) $$|\langle \widetilde{u}, \ker (\theta_{{i_0}}I_n-P_n) \rangle|^{2} \convas  {\f{-1}{\theta_{i_0}^2G_{\mu_X}'(\rho)}}$$
where $\rho$ is the limit of $\widetilde{\lambda}_{\lceil nc_{j-1} \rceil +i_0}$ given by Theorem \ref{140709.main};

\flushleft (b)
$$\langle \widetilde{u}, \oplus_{i\neq {i_0}} \ker (\theta_{i}I_n-P_n) \rangle \convas 0.$$
 \end{Th}
\begin{proof}
The result is obtained by following the approach taken in \cite[pp. 514-516]{benaych2011eigenvalues} for proving Theorem 2.2 and accounting for the possibly multi-valued nature of $G_{\mu_X}^{-1}(1/\theta)$.
\end{proof}

\begin{Th}[Phase transition of eigenvector informativeness]\label{200709.13h5}When $r=1$, let the sole non-zero eigenvalue of $P_n$  be denoted by $\theta$.
Consider $j \in \{1, \ldots, \ell \}$ such that
$$ \ff{\theta}\notin ( G_{\mu_X}(a_{j-1}^-),G_{\mu_X}(b_j^+)),\quad\textrm{and}\quad G_{\mu_X}'(b_{j}^+)=-\infty \textrm{ and } G_{\mu_X}(a_{j-1}^-)=-\infty.$$
For each $n$, let $\widetilde{u}$ be a unit-norm eigenvector of $\wtX$ associated with $\widetilde{\lambda}_{\lceil nc_{j-1} \rceil +1}$. Then we have $$\lan \widetilde{u}, \ker(\theta I_n-P_n)\ran\convas 0, $$
as $n \longrightarrow \infty$.
\end{Th}
\begin{proof}
The result is obtained by following the approach taken in \cite[pp. 516-517]{benaych2011eigenvalues} for proving Theorem 2.3 and accouting for the possibly multi-valued nature of $G_{\mu_X}^{-1}(1/\theta)$.
\end{proof}

The following proposition allows to assert that in many classical matrix models, such as Wigner or Wishart matrices, the above phase transitions actually occur with a finite threshold.

\begin{propo}[Edge density decay condition and the phase transition]\label{square root add}
Assume that the limiting eigenvalue distribution $\mu_X$, supported on $\ell$ disjoint intervals,  has a density $f_{\mu_X}$ with a power decay at $b_j$ for $j = 1, \ldots, \ell$, i.e., that, as $t\to b_j$ with $t<b_j$,  $f_{\mu_X}(t) \sim c (b_j-t)^\alpha$  for some exponent $\alpha>-1$ and  some constant $c$. Then:
$$ G_{\mu_X}(b_j^{+}) < \infty\iff \alpha>0  \qquad\textrm{ and } \qquad G'_{\mu_X}(b_j^{+})=-\infty \iff \al\le 1.$$
Similarly, if $f_{\mu_X}$ has a power decay at $a_{j-1}$ for $j = 2, \ldots, \ell$, i.e., that, as $t\to a_{j-1}$ with $t>a_{j-1}$,  $f_{\mu_X}(t) \sim c (t-a_{j-1})^\alpha$  for some exponent $\alpha>-1$ and  some constant $c$. Then
$$ G_{\mu_X}(a_{j-1}^{-}) < \infty\iff \alpha>0  \qquad\textrm{ and } \qquad G'_{\mu_X}(a_{j-1}^{-})=-\infty \iff \al\le 1.$$
\end{propo}

Theorem \ref{140709.main} describes the fundamental limits of eigen-gap based signal detection. Principal eigen-gap detection will fail whenever $\theta_i < 1/G_{\mu_X}(b_{1}^+)$. If  $G_{\mu_X}(b_{j}^+) > G_{\mu_X}(b_{1}^+)$ for $j = 2, \ldots, \ell$ then principal eigen-gap detection will be suboptimal as the middle eigen-gaps will reveal the presence of a low-rank signal even when the principal eigen-gap does not. Theorem \ref{180709.13h39} shows that whenever there is an eigen-gap, the corresponding eigenvectors will be informative. Theorem \ref{200709.13h5} provides insight on the fundamental limits of low-rank signal matrix estimation.

\subsection{Singular values and singular vectors}\label{sec:main sv}

Let  $X_n$ be  an $n \times m$ ($n \leq m$, without loss of generality) random matrix whose ordered singular values we denote by $\si_{1}(X_n) \geq \cdots \geq \si_{n}(X_n)$. Let $\mu_{X_{n}}$ be the empirical singular value distribution, \ie, the probability measure defined as
\[
\mu_{X_{n}} = \frac{1}{n} \sum_{i=1}^{n} \delta_{\si_{i}(X_{n})}.
\]
As before, assume that the probability measure $\mu_{X_{n}}$ converges almost surely weakly, as $n \longrightarrow \infty$, to a non-random compactly supported  \pro measure $\mu_X$ that is supported on $\ell$ disjoint intervals so that
$$\mu_X = \sum_{j=1}^{\ell} p_{j}\, \mu_{j} ,$$
where for $j = 1, \ldots, \ell$, the measures $\mu_{j}(x)$ is a non-random probability measure are supported on $[a_{j},b_{j}]$ with  $d\mu_{j}(z) > 0$ for $z \in (a_{j},b_{j})$ and $a_{\ell} < b_{\ell} < a_{\ell - 1} < b_{\ell-1} < \ldots  < a_{1} < b_{1}$. Define $c_{j} = \sum_{i=0}^{j} p_{i}$ with $p_0 := 0$, $c_0 := 0$ and $c_{\ell} := 1$. We assume that  for $j = 1, \ldots,  \ell$, $\si_{\lceil n c_{j-1} \rceil + 1} \convas b_{j}$ and that $\si_{\lceil  n c_{j} \rceil} \convas a_{i}$. As before, we use $\lceil n c_{j-1} \rceil$ to denote the smallest integer greater than or equal to $n c_{j-1}$.

For a given $r \ge 1$, let $\tta_1\ge\cdots\ge \tta_r$ be deterministic non-zero real numbers, chosen independently of $n$.  For every $n$, let $P_n$ be an $n \times m$ matrix   having rank $r$ with its $r$ non-zero singular values equal to $\theta_{1}, \ldots, \theta_{r}$. We suppose that $X_n$ and $P_n$ are independent and that $X_n$, the noise-only matrix  is bi-unitarily invariant while the low-rank signal matrix $P_n$ is deterministic.  Recall that a random matrix is said to be {\it bi-orthogonally invariant} (or {\it bi-unitarily invariant}) if its distribution is invariant under multiplication on the left and right by orthogonal (or unitary) matrices. Alternately, if $P_n$ has isotropically random right (or left) singular vectors then, then $X_n$ need not be unitarity invariant under multiplication on the right (or left, respc.) by orthogonal or unitary matrices. Equivalently, $X_n$ can have deterministic right and left singular vectors while $P_n$ can have isotropically random left and right singular vectors and we would get the same result stated shortly.

Consider the rank $r$ additive perturbation of the random matrix $X_n$ given by
$$\wtX_n = X_n+P_n.$$
where
$$P_n = \sum_{i=1}^{r} \theta_{i} u_{i}  v_{i}^{*},$$
and $\{ u_{i} \}_{i=1}^{r}$ and $\{ v_{i} \}_{i=1}^{r}$ are the left and right singular vectors, respectively of $P_n$.

For this model, we establish the following results.

\begin{Th}[Largest singular value phase transition]\label{140709.main.rectangular}
The singular values of $\wtX $ exhibit the following    behavior as $n,m_n \to \infty$ and $n/m_n \to c$. . We have that for each $1\le  i \leq r$ and $ 1 \leq j \leq \ell$,
$$\si_{\lceil nc_{j-1} \rceil +i}(\wtX) \convas \begin{cases} D_{\mu_X,(b_j,a_j-1)}^{-1}(1/\theta_i^2)& 1/ D_{\mu_X}(b_j^+)  <\theta_i^2 < 1/ D_{\mu_X}(a_{j-1}^-) ,\\ \\
b_j &\textrm{ if } \theta_{i}^2 < 1/D_{\mu_X}(b_j^+), \\ \\
a_{j-1} & \textrm{ if } \theta_{i}^2 > 1/D_{\mu_X}(a_{j-1}^-),
\end{cases}$$
where $D_{\mu_X}$, the $D$-transform of $\mu_X$ defined by
$$ D_{\mu_X}(z) :=\left[\int \frac{z}{z^2-t^{2}} \ud \mu_X(t)\right] \times \left[c\int \frac{z}{z^2-t^{2}} \ud \mu_X(t)+\frac{1-c}{z}\right] \qquad \textrm{for } z \notin \cup_{j=1}^{\ell} [a_j,b_{j}],$$
and   $D_{\mu_X,(b_j,a_{j-1})}^{-1}(\cdot)$ will denote its functional inverse on $(b_j,a_{j-1})$ with $a_{0}  = + \infty$.
\end{Th}
\begin{proof}
The result is obtained by following the approach taken in \cite[pp. 127--129]{benaych2012singular} for proving Theorem 2.9 and accounting for the possibly multi-valued nature of the $D_{\mu_X}^{-1}(\cdot)$. The  key ingredient of the proof is the recognition that the non-zero, positive eigenvalues of 
$$ 
\begin{bmatrix} 0 & \wtX \\ \wtX^* & 0 \end{bmatrix} = 
\begin{bmatrix} 0 & X \\ X^* & 0 \end{bmatrix}+
\begin{bmatrix} 0 & \sum_{i=1}^{r} \theta_{i} u_i v_i^*\\ \sum_{i=1}^{r} \theta_{i} v_i u_i^* & 0 \end{bmatrix},$$
are precisely the singular values of $\wtX$. Thus adopting the approach outlined in Section \ref{sec:master equations} while taking into account the structured rank $2r$ perturbation gives us the stated result.
\end{proof}

\begin{Th}[Informativeness of singular vectors]\label{180709.13h39.rectangular}
 Assume throughout that $\theta >0$ and let $D_{\mu_X}(a_0^-) = +\infty$. Consider   ${i_0}\in \{1, \ldots, r\}$  \st $1/\theta_{i_0}^2 \in \bigcup_{j = 1}^{l} (D_{\mu_X}(a_{j-1}^-),D_{\mu_X}(b_{j}^+), )$. For each such $i_0$, consider $j(i_0) \equiv j \in \{1, \ldots, l\}$ such that $1/\theta_{i_0}^2 \in (D_{\mu_X}(a_{j-1}^{-}),D_{\mu_X}(b_j^+))$ and let $\widetilde{u}$  and $\widetilde{v}$ be unit-norm left and right singular vectors of $\wtX$  associated with the singular value $\widetilde{\sigma}_{\lceil nc_{j-1} \rceil+i_0} $. Then we have, as $n \longrightarrow \infty$,

\flushleft a)
\be\label{250709.09h13}
|\langle \widetilde{u}, \Span\{u_i\ste \tta_i=\tta_{i_0}\}\rangle|^{2}  \convas {\f{-2\vfi_{\mu_X}(\rho)}{\theta_{i_0}^2 D'_{\mu_X}(\rho) }},
\ee
\flushleft b)
\be\label{250709.09h13.300709}
|\langle \widetilde{v}, \Span\{v_i\ste \tta_i=\tta_{i_0}\}  \rangle|^{2}  \convas  {\f{-2\vfi_{\widetilde{\mu}_X}(\rho)}{\theta_{i_0}^2 D_{\mu_X}'(\rho) }},
\ee
 where $\rho$ is the limit of $\widetilde{\si}_{i_0}$ given by Theorem \ref{140709.main.rectangular} and $\widetilde{\mu}_X=c\mu_X+(1-c)\delta_0$ and for any \pro measure $\mu$,
\be\label{9710.19h07}\vfi_\mu(z):=\int \f{z}{z^2-t^2}\ud \mu(t).\ee
\flushleft c) Furthermore, in the same asymptotic limit, we have
$$ |\langle \widetilde{u}, \Span\{u_i\ste \tta_i\ne \tta_{i_0}\} \rangle|^{2}    \convas 0, \quad \textrm{ and } \quad |\langle \widetilde{v},  \Span\{v_i\ste \tta_i\ne \tta_{i_0}\} \rangle|^{2} \convas 0 ,$$
and
$$\lan \vfi_{{\mu}_X}(\rho){P_n}\widetilde{v}-\widetilde{u} \;,\;  \Span\{u_i\ste \tta_i=\tta_{i_0}\} \rangle \convas 0.$$
\end{Th}
\begin{proof}
The result is obtained by following the approach taken in \cite[pp. 129--131]{benaych2012singular} for proving Theorem 2.10 and accounting for the possibly multi-valued nature of the $D_{\mu_X}^{-1}(\cdot)$.
\end{proof}

\begin{Th}[Phase transition of vector informativeness]\label{deloc.sing.vect.250709.rectangular}When $r=1$, let the  sole singular value of $P_n$  be denoted by $\theta$.
Consider $j \in \{1, \ldots, \ell \}$ such that
$$ \ff{\theta^2} \notin ( D_{\mu_X}(a_{j-1}^-),D_{\mu_X}(b_j^+)),\quad\textrm{and}\quad D_{\mu_X}'(b_{j}^+)=-\infty \textrm{ and  } D_{\mu_X}'(a_{j-1}^-)=-\infty.$$
For each $n$, let $\widetilde{u}$ and $\widetilde{v}$ be  unit-norm left and right singular vectors of $\wtX$ associated with $\widetilde{\sigma}_{ \lceil nc_{j-1} \rceil+1}$. Then we have that
$$ \langle \widetilde{u},   \ker (\theta^2 I_n-P_nP_n^*)\rangle \convas 0, \quad \textrm{ and } \quad\langle \widetilde{v},   \ker (\theta^2I_m-P_n^*P_n)\rangle \convas 0,$$
as $n \longrightarrow \infty$.
\end{Th}
\begin{proof}
The result is obtained by following the approach taken in \cite[pp. 131]{benaych2012singular} for proving Theorem 2.11 and accounting for the possibly multi-valued nature of the $D_{\mu_X}^{-1}(\cdot)$.
\end{proof}

Theorem \ref{140709.main} describes the fundamental limits of eigen-gap based signal detection. Principal gap detection will fail whenever $\theta_i^2 < 1/D_{\mu_X}(b_{1}^+)$. If  $D_{\mu_X}(b_{j}^+) > D_{\mu_X}(b_{1}^+)$ for $j = 2, \ldots, \ell$ then principal eigen-gap detection will be suboptimal as the middle eigen-gaps will reveal the presence of a low-rank signal even when the principal eigen-gap does not. Theorem \ref{180709.13h39.rectangular} shows that whenever there is an eigen-gap, the corresponding singular vectors will be informative. Theorem \ref{deloc.sing.vect.250709.rectangular} provides insight on the fundamental limits of low-rank signal matrix estimation. The analog of Proposition \ref{square root add} also applies here.

\section{Noise models that might produce informative middle components}\label{sec:noise models}

Our discussion has brought into sharp focus the pivotal role played by the noise eigen-spectrum in determining the relative informativeness of the principal and middle components  of the singular value (or eigen) decomposition of signal-plus-noise data matrix models as in (\ref{eq:sigplusnoise mat}).

Specifically, we showed that if the noise eigen-spectrum is supported on a single connected interval then the principal components will indeed (with high probability) be the most informative components and their use in detection and estimation is justified.

However, if the noise eigen-spectrum is supported on multiple intervals, as in Figure \ref{fig:mca phase explanation}, then the  principal components will remain informative in the high SNR regime (\textit{i.e.}, large $\theta_i$). However, for moderate to low SNR, the middle components might also be informative and may remain informative even when the principal components are no longer informative. In such settings, identifying large middle eigen-gaps and using the associated middle eigenvectors for inference can improve inference .

This leads to a natural question: When will the noise eigen-spectrum exhibit a disconnected spectrum?

We conclude by identifying a large class of Gaussian mixture models that produce precisely such an eigen-spectrum.  Consider the class of noise matrices modeled as
$$X = G \Sigma^{1/2},$$
where $G$ is an $m \times n$ matrix with i.i.d. mean zero, variance $1/m$ (say) Gaussian entries. If the rows of $X$ denote spatial measurements and the columns represent temporal measurements, then $\Sigma$ is a temporal covariance matrix and $XX^*$ is a Wishart distributed matrix.  These models arise in many statistical signal processing and machine learning applications where PCA/SVD is often used as the first step in inferential process (see, for e.g. \cite{tipping1999mixtures,dasgupta1999learning,sanjeev2001learning,vempala2002spectral,kannan2005spectral,hsu2012learning}).

\begin{figure}[t]
\vspace{-0.15cm}
\centering
\includegraphics[trim = 0 170 0 170, clip = true,width=5.85in]{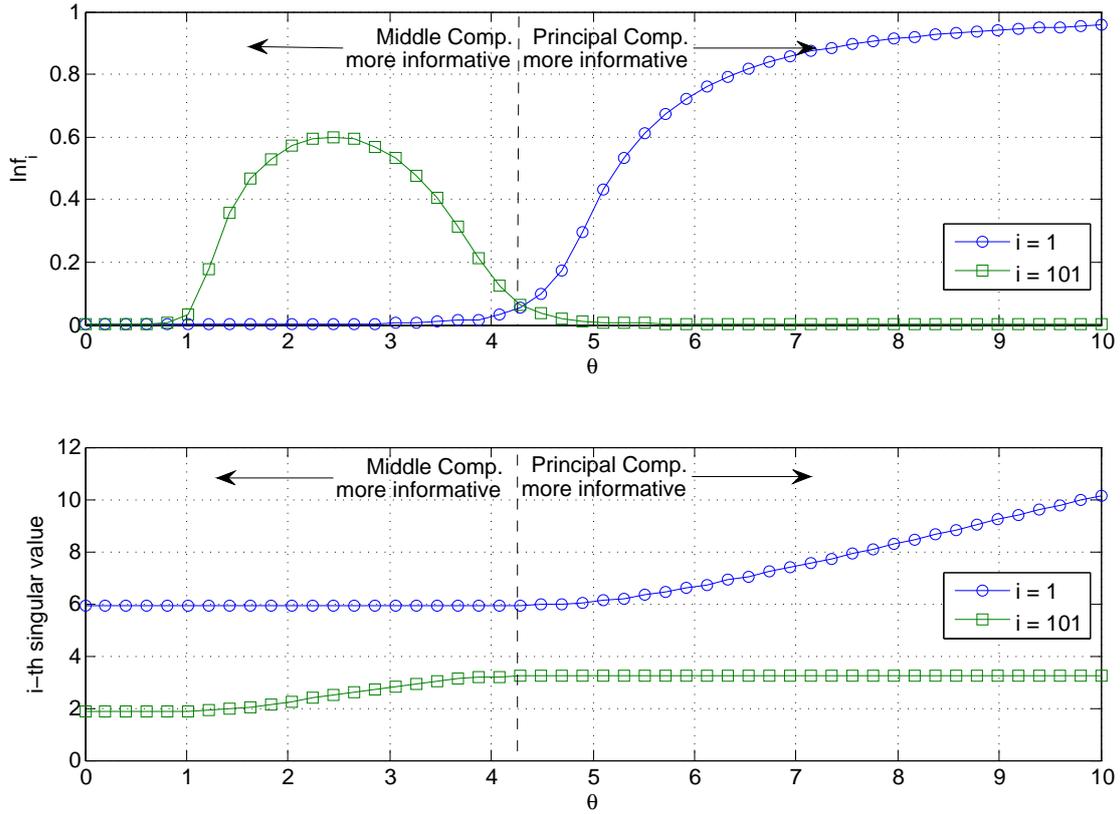}
\caption{The evolution of the informativeness of the principal and the middle component as a function of $\theta$ for the model in (\ref{eq:sigplusnoise mat}) and $X =  G \Sigma^{1/2}$, where $G$ in an $n \times m$ matrix (here $n = m = 1000$) with i.i.d. mean zero, variance $1/m$ entries and $\Sigma = \textrm{diag}(20 I_{n/10}, I_{n - n/10})$. The upper panel plots  $\inform{i} = |\langle \widetilde{u}_{i},u \rangle|^2$ computed over $250$ Monte-Carlo trials. The lower panel plots the $i$-th largest singular value of $\wtX$. }
\label{fig:mca numerical expt}
\end{figure}

Bai and Silverstein characterize the limiting eigenvalue distribution of $XX^*$ in \cite{silverstein1995empirical}. What emerges from their analysis \cite{bai1998no,silverstein1992signal} is that for the noise eigen-spectrum of $XX^*$ to have a disconnected spectrum, the eigenvalue spectrum of $\Sigma$ has to have a limiting distribution that is supported on disconnected intervals. In addition, the separation between these intervals has to be relatively large for the spectrum of $X$ to be supported on disconnected intervals. There is no simple formula for how large this separation has to be. There are expressions in \cite{silverstein1995empirical} for the form of the spectrum of $X$ as a function of the spectrum $\Sigma$, from which it can be ascertained whether the support is supported on multiple intervals on not using the results in \cite{bai1998no,silverstein1992signal}. Moreover, the spectrum will exhibit square-root decay at the edges \cite{silverstein1995analysis} and so the phase transitions described will manifest.

 Figure \ref{fig:mca numerical expt} plots evolution of the $i$-th singular value and $\inform{i} = |\langle \widetilde{u}_{i},u \rangle|^2$ as a function of $\theta$ for the model in (\ref{eq:sigplusnoise mat}) with $X = G \Sigma^{1/2}$ and $\Sigma = \textrm{diag}(20 I_{n/10}, I_{n - n/10})$. The figure clearly shows the phase transition in the informativness of the principal and middle components and shows that there is a low SNR regime where the middle component is informative even when the principal component is not.  The values where the phase transitions occur can be theoretically predicted, if so desired, using Theorem \ref{140709.main.rectangular}. Figure \ref{fig:mca sv} shows a sample realization of the singular values for the same setting when $\theta = 2$.

Thus the potential for informative middle components to emerge is the greatest in large, heterogeneous datasets where there might be significant temporal (or spatial) variation. These might be exploited for extracting additional processing gain beyond what principal component analysis might offer.

Conversely, if the temporal covariance matrix $\Sigma$ represents a relatively homogenous (in time) data set, then there will be no gap in the eigen-spectrum and the use of principal components is justifiably optimal.

Expanding the range of noise models for which similar predictions can be made is a natural next step.  It remains an open problem to fully characterize the vanishing informativeness of the components of the singular value decomposition associated with singular values that (asymptotically) exhibit an $o(1)$ eigen-gap.  Additional hypotheses on the noise eigen-distribution will likely be required - establishing natural conditions for these remains an important line of inquiry. A result along these lines would firmly establish that the informative components associated with the singular/eigen values that exhibit an eigen-gap are indeed the \textit{maximally} informative components.

\clearpage

\bibliographystyle{plain}
\bibliography{spiked_bib}

\begin{thebibliography}{10}

\bibitem{bai1998no}
ZD~Bai and J.W. Silverstein.
\newblock No eigenvalues outside the support of the limiting spectral
  distribution of large-dimensional sample covariance matrices.
\newblock {\em The Annals of Probability}, 26(1):316--345, 1998.

\bibitem{baik2005phase}
J.~Baik, G.~Ben~Arous, and S.~P{\'e}ch{\'e}.
\newblock Phase transition of the largest eigenvalue for nonnull complex sample
  covariance matrices.
\newblock {\em The Annals of Probability}, 33(5):1643--1697, 2005.

\bibitem{baik2006eigenvalues}
J.~Baik and J.W. Silverstein.
\newblock Eigenvalues of large sample covariance matrices of spiked population
  models.
\newblock {\em Journal of Multivariate Analysis}, 97(6):1382--1408, 2006.

\bibitem{benaych2011eigenvalues}
F.~Benaych-Georges and R.R. Nadakuditi.
\newblock The eigenvalues and eigenvectors of finite, low rank perturbations of
  large random matrices.
\newblock {\em Advances in Mathematics}, 227(1):494--521, 2011.

\bibitem{benaych2012singular}
F.~Benaych-Georges and R.R. Nadakuditi.
\newblock The singular values and vectors of low rank perturbations of large
  rectangular random matrices.
\newblock {\em Journal of Multivariate Analysis}, pages 120--135, 2012.

\bibitem{brand2006fast}
M.~Brand.
\newblock Fast low-rank modifications of the thin singular value decomposition.
\newblock {\em Linear algebra and its applications}, 415(1):20--30, 2006.

\bibitem{dasgupta1999learning}
S.~Dasgupta.
\newblock Learning mixtures of gaussians.
\newblock In {\em Foundations of Computer Science, 1999. 40th Annual Symposium
  on}, pages 634--644. IEEE, 1999.

\bibitem{deshpande2006adaptive}
A.~Deshpande and S.~Vempala.
\newblock Adaptive sampling and fast low-rank matrix approximation.
\newblock {\em Approximation, Randomization, and Combinatorial Optimization.
  Algorithms and Techniques}, pages 292--303, 2006.

\bibitem{drineas2005nystrom}
P.~Drineas and M.W. Mahoney.
\newblock On the nystr{\"o}m method for approximating a gram matrix for
  improved kernel-based learning.
\newblock {\em The Journal of Machine Learning Research}, 6:2153--2175, 2005.

\bibitem{eckart1936approximation}
C.~Eckart and G.~Young.
\newblock The approximation of one matrix by another of lower rank.
\newblock {\em Psychometrika}, 1(3):211--218, 1936.

\bibitem{el2007tracy}
N.~El~Karoui.
\newblock Tracy--widom limit for the largest eigenvalue of a large class of
  complex sample covariance matrices.
\newblock {\em The Annals of Probability}, 35(2):663--714, 2007.

\bibitem{friedman2001elements}
J.~Friedman, T.~Hastie, and R.~Tibshirani.
\newblock {\em The elements of statistical learning}, volume~1.
\newblock Springer Series in Statistics, 2001.

\bibitem{frieze2004fast}
A.~Frieze, R.~Kannan, and S.~Vempala.
\newblock Fast monte-carlo algorithms for finding low-rank approximations.
\newblock {\em Journal of the ACM (JACM)}, 51(6):1025--1041, 2004.

\bibitem{halko2011finding}
N.~Halko, P.G. Martinsson, and J.A. Tropp.
\newblock Finding structure with randomness: Probabilistic algorithms for
  constructing approximate matrix decompositions.
\newblock {\em SIAM review}, 53(2):217--288, 2011.

\bibitem{hp00}
Fumio Hiai and D{\'e}nes Petz.
\newblock {\em The semicircle law, free random variables and entropy},
  volume~77 of {\em Mathematical Surveys and Monographs}.
\newblock American Mathematical Society, Providence, RI, 2000.

\bibitem{hsu2012learning}
D.~Hsu and S.M. Kakade.
\newblock Learning gaussian mixture models: Moment methods and spectral
  decompositions.
\newblock {\em arXiv preprint arXiv:1206.5766}, 2012.

\bibitem{johnstone2001distribution}
I.M. Johnstone.
\newblock On the distribution of the largest eigenvalue in principal components
  analysis.(english.
\newblock {\em Ann. Statist}, 29(2):295--327, 2001.

\bibitem{johnstone2006high}
I.M. Johnstone.
\newblock High dimensional statistical inference and random matrices.
\newblock In {\em Proceedings oh the International Congress of Mathematicians:
  Madrid, August 22-30, 2006: invited lectures}, pages 307--333, 2006.

\bibitem{jolliffe2005principal}
I.~Jolliffe.
\newblock {\em Principal component analysis}.
\newblock Wiley Online Library, 2005.

\bibitem{kannan2005spectral}
R.~Kannan, H.~Salmasian, and S.~Vempala.
\newblock The spectral method for general mixture models.
\newblock {\em Learning Theory}, pages 155--199, 2005.

\bibitem{kritchman2008determining}
S.~Kritchman and B.~Nadler.
\newblock Determining the number of components in a factor model from limited
  noisy data.
\newblock {\em Chemometrics and Intelligent Laboratory Systems}, 94(1):19--32,
  2008.

\bibitem{kritchman2009non}
S.~Kritchman and B.~Nadler.
\newblock Non-parametric detection of the number of signals: hypothesis testing
  and random matrix theory.
\newblock {\em Signal Processing, IEEE Transactions on}, 57(10):3930--3941,
  2009.

\bibitem{mirsky1960symmetric}
L.~Mirsky.
\newblock Symmetric gauge functions and unitarily invariant norms.
\newblock {\em The quarterly journal of mathematics}, 11(1):50--59, 1960.

\bibitem{nadakuditi2008sample}
R.R. Nadakuditi and A.~Edelman.
\newblock Sample eigenvalue based detection of high-dimensional signals in
  white noise using relatively few samples.
\newblock {\em Signal Processing, IEEE Transactions on}, 56(7):2625--2638,
  2008.

\bibitem{nadler2010nonparametric}
B.~Nadler.
\newblock Nonparametric detection of signals by information theoretic criteria:
  performance analysis and an improved estimator.
\newblock {\em Signal Processing, IEEE Transactions on}, 58(5):2746--2756,
  2010.

\bibitem{onatski2010determining}
A.~Onatski.
\newblock Determining the number of factors from empirical distribution of
  eigenvalues.
\newblock {\em The Review of Economics and Statistics}, 92(4):1004--1016, 2010.

\bibitem{paul2007asymptotics}
D.~Paul.
\newblock Asymptotics of sample eigenstructure for a large dimensional spiked
  covariance model.
\newblock {\em Statistica Sinica}, 17(4):1617, 2007.

\bibitem{sanjeev2001learning}
A.~Sanjeev and R.~Kannan.
\newblock Learning mixtures of arbitrary gaussians.
\newblock In {\em Proceedings of the thirty-third annual ACM symposium on
  Theory of computing}, pages 247--257. ACM, 2001.

\bibitem{silverstein1995empirical}
J.W. Silverstein and ZD~Bai.
\newblock On the empirical distribution of eigenvalues of a class of large
  dimensional random matrices.
\newblock {\em Journal of Multivariate analysis}, 54(2):175--192, 1995.

\bibitem{silverstein1995analysis}
J.W. Silverstein and S.I. Choi.
\newblock Analysis of the limiting spectral distribution of large dimensional
  random matrices.
\newblock {\em Journal of Multivariate Analysis}, 54(2):295--309, 1995.

\bibitem{silverstein1992signal}
J.W. Silverstein and P.L. Combettes.
\newblock Signal detection via spectral theory of large dimensional random
  matrices.
\newblock {\em Signal Processing, IEEE Transactions on}, 40(8):2100--2105,
  1992.

\bibitem{tipping1999mixtures}
M.E. Tipping and C.M. Bishop.
\newblock Mixtures of probabilistic principal component analyzers.
\newblock {\em Neural computation}, 11(2):443--482, 1999.

\bibitem{vempala2002spectral}
S.~Vempala and G.~Wang.
\newblock A spectral algorithm for learning mixtures of distributions.
\newblock In {\em Foundations of Computer Science, 2002. Proceedings. The 43rd
  Annual IEEE Symposium on}, pages 113--122. IEEE, 2002.

\end{thebibliography}

\end{document}